\documentclass{article}
\usepackage{latexsym,amsmath,amssymb,amscd}
\usepackage[dvips]{color}
\numberwithin{equation}{section}

\begin{document}

\title{A Note on the Axisymmetric Stationary Metric in the General Theory of Relativity}
\author{Tetu Makino \footnote{Professor Emeritus at Yamaguchi University, Japan. E-mail: makino@yamaguchi-u.ac.jp}}
\date{\today}
\maketitle

\newtheorem{Lemma}{Lemma}
\newtheorem{Proposition}{Proposition}
\newtheorem{Theorem}{Theorem}
\newtheorem{Definition}{Definition}
\newtheorem{Remark}{Remark}
\newtheorem{Corollary}{Corollary}

\begin{abstract}
We consider the equations for the coefficients of stationary rotating axisymmetric metrics governed by the Einstein-Euler equations, that is, the Einstein equations together with the energy-momentum tensor of a barotropic perfect fluid. 
Although the derived equations are not already known except for the case of the constant angular velocity described in the corotating coordinate system, the main content of this article is not to derive the equations, but to prove the equivalence of the derived equations with the full set of the Einstein equations, and to prove the consistency of the derived equations. These affairs have not yet been discussed except for the vacuum case, and are far from being self-evident, requiring tedious careful calculations and some tricks. The proof is done under the assumption that the  angular velocity is constant on a neighborhood of the support of the density. The conclusions seem to be doubtful if this assumption does not hold.\\


\end{abstract}

\section{Introduction}

We consider the metric 
\begin{equation}
ds^2=g_{\mu\nu}dx^{\mu}dx^{\nu}
\end{equation}
satisfying the Einstein equations
\begin{equation}
R_{\mu\nu}-\frac{1}{2}g_{\mu\nu}R=
\frac{8\pi\mathsf{G}}{\mathsf{c}^4}T_{\mu\nu}.
\end{equation}
The positive constants $\mathsf{c}, \mathsf{G}$ are the speed of light and the constant of gravitation. 
The energy-momentum tensor $T^{\mu\nu}$ is supposed to be that of a perfect fluid
\begin{equation}
T^{\mu\nu}=(\epsilon +P)U^{\mu}U^{\nu}-Pg^{\mu\nu}
\quad
\mbox{with}\quad \epsilon =\mathsf{c}^2\rho.
\end{equation}

Here we suppose\\

{\bf (A)}: {\it The pressure $P$ is  a  smooth function of $\rho>0$ such that
 $0<P,
0<dP/d\rho <\mathsf{c}^2$ for $\rho >0$ and there are constants $\gamma, \mathsf{A}$ and a function $\Upsilon$ which is analytic near $0$ and satisfies $\Upsilon(0)=0$ such that
\begin{equation}
P=\mathsf{A}\rho^{\gamma}(1+\Upsilon(\mathsf{A}\rho^{\gamma-1}/\mathsf{c}^2)) 
\end{equation}
for $\rho > 0$ and $1<\gamma < 2, 0<\mathsf{A}$.
 }\\

Also, we consider the dust:
$P=0$ for $\forall \rho$, too.\\

In this Note we consider axially symmetric metrics, that is, taking
the coordinates
\begin{equation}
x^0=\mathsf{c}t,\quad x^1=\varpi,\quad x^2=\phi,\quad x^3=z,
\end{equation}
we consider the metric $ds^2=g_{\mu\nu}dx^{\mu}dx^{\nu}$ in the following form (Lewis 1932 \cite{Lewis}, 
Papapetrou 1966 \cite{Papapetrou} ), say, `{\bf Lanczos form}' after \cite{Lanczos}: 
\begin{equation}
ds^2=e^{2F}
(\mathsf{c}dt+Ad\phi)^2-
e^{-2F}[e^{2K}
(d\varpi^2+dz^2)+\Pi^2d\phi^2], \label{2}
\end{equation}
that is,
\begin{align*}
& g_{00}=e^{2F}, \quad g_{02}=g_{20}=e^{2F}A, \quad
g_{11}=g_{33}=-e^{2(-F+K)}, \\
&g_{22}=e^{2F}A^2-e^{-2F}\Pi^2,\quad \mbox{other}\  g_{\mu\nu}=0,
\end{align*}
where the quantities $F, A, K$ and $\Pi$ are functions of only  $\varpi$ and $z$.
$\mathsf{c}^2F$ is called the `{\bf generalized Newtonian potential}' and $\mathsf{c}A$ is called
`{\bf gravitomagnetic potential}'.

The components of the metric and the Christoffel symbols 
$$ \Gamma_{\nu\lambda}^{\mu}=\frac{1}{2}g^{\mu\alpha}
(\partial_{\lambda}g_{\alpha\nu}+
\partial_{\nu}g_{\alpha\lambda}-\partial_{\alpha}g_{\nu\lambda} )$$
are given in {\bf Appendix 1}.
Here and hereafter $\partial_{\mu}$ stands for 
$\displaystyle \frac{\partial}{\partial x^{\mu}}$, while $\nabla_{\mu}$ will stand for
the covariant derivative
$$\nabla_{\lambda}A^{\mu\nu}=
\partial_{\lambda}A^{\mu\nu}
+\Gamma_{\alpha\lambda}^{\mu}A^{\alpha\nu}+
\Gamma_{\alpha\lambda}^{\nu}A^{\mu\alpha}.
$$\\

The 4-velocity vector filed $U^{\mu}$ is supposed to be of the form
\begin{equation}
U^{\mu}\frac{\partial}{\partial x^{\mu}}=e^{-G}\Big(\frac{1}{\mathsf{c}}\frac{\partial}{\partial t}+
\frac{\Omega}{\mathsf{c}}\frac{\partial}{\partial\phi}\Big),
\end{equation}
that is,
\begin{equation}
U^0=e^{-G},\quad U^1=U^3=0,\quad
U^2=e^{-G}\frac{\Omega}{\mathsf{c}},
\end{equation}
where $\Omega$ depends only on $\varpi, z$.

Since $U^{\mu}U_{\mu}=1$, the factor $e^{-G}=U^0$ is given by
\begin{equation}
e^{2F}\Big(1+\frac{\Omega}{\mathsf{c}}A\Big)^2
-e^{-2F}\frac{\Omega^2}{\mathsf{c}^2}\Pi^2=e^{2G}. \label{18}
\end{equation}

In this point of view, $F, A, \Pi$ to be considered should satisfy the following assumption 
{\bf (B1) } with respect to the prescribed $\Omega$ :\\

{\bf (B1)}: {\it It holds that}
\begin{equation}
e^{2F}\Big(1+\frac{\Omega}{\mathsf{c}}A\Big)^2-e^{-2F}
\frac{\Omega^2}{\mathsf{c}^2}\Pi^2 >0.
\end{equation}\\

The plan of this note is as follows. In Section 2 the full set of Einstein equations is described by the components $F, A, K, \Pi$ of the metric \eqref{2}, and by the auxiliary quantities
$f, k, m, l$, which are useful alternatives of $F, A, K, \Pi$. In Section 3  the equations for
$F, A, K, \Pi$ are described. The process of derivation is clarified. The equivalence of this reduced system of equations and the original full set of Einstein equations is discussed. Moreover the consistency of the first order system for $K$ is discussed. To these problems the affirmative answers are given under the restriction that either $\Omega$ is constant or $\rho=0$ on the considered domain. The self-contained proofs of these results are new. In Section 4 the corresponding discussion is described concerning the associated quantities $F', A', K', \Pi'(=\Pi) $. These quantities can be interpreted as the components of the metric alternative to $F, A, K, \Pi$ with respect to the corotating coordinate system when $\Omega$ is constant, as noted in Section 5. But here we consider variable $\Omega$, say, the case of differentially rotating stars. Although a part of results for the case of constant $\Omega$ was already given in the preceding work \cite{asEE}, the most part is new. Actually, during the preceding study of \cite{asEE}, the author lacked the sense of necessity to prove the equivalence between the subset of 3 equations in the full Einstein equations concerning $(\mu,\nu)=(1, 1), (3,3), (1,3) $ and the set of 2 equations with respect to $K'$. Later the author was aware that the equivalence is announced in \cite[p. 27]{Islam} for the case of vacuum, say,  when $\rho=P=0$ and $\Pi=\varpi$. Anyway we should note that the problems of the equivalence and the consistency require tedious calculations with some tricks, and are far from being self-evident. This work is motivated by the matter-vacuum matching problem, on which the author gave a report and received helpful comments in the  BIRS-CMO Workshop `Time-like Boundaries in General Relativistic Evolution Problems' , July 28 - August 2, 2019.

\section{Einstein equations}

\subsection{Euler equations and the variable $u$}

The non-zero components of the energy-momentum tensor $T^{\mu\nu}$ are:
\begin{align*}
&T^{00}=e^{-2G}(\epsilon+P)+
\frac{P}{\Pi^2}(e^{2F}A^2-e^{-2F}\Pi^2), \\
&T^{02}=T^{20}=
e^{-2G}(\epsilon+P)\frac{\Omega}{\mathsf{c}}-
\frac{P}{\Pi^2}e^{2F}A, \\
&T^{11}=T^{33}=Pe^{2F-2K}, \\
&T^{22}=
e^{-2G}(\epsilon+P)\frac{\Omega^2}{\mathsf{c}^2}+
\frac{P}{\Pi^2}e^{2F}.
\end{align*}

The Euler equations are $\nabla_{\mu}T^{\mu\nu}=0$. 
Let us consider them.

Keeping in mind \eqref{18}, we can show through tedious calculations that 
$$\nabla_{\mu}T^{\mu 1}=
e^{2F-2K}\Big(\partial_1P+(\epsilon+P)\partial_1G
-(\epsilon+P)U^0U_2\frac{\partial_1\Omega}{\mathsf{c}}\Big) $$
and
$$\nabla_{\mu}T^{\mu 3}=
e^{2F-2K}\Big(\partial_3P+(\epsilon+P)\partial_3G
-(\epsilon+P)U^0U_2\frac{\partial_3\Omega}{\mathsf{c}}\Big). $$
Here
\begin{align}
U^0U_2&=e^{2F-2G}\Big(A\Big(1+\frac{\Omega}{\mathsf{c}}\Big)-e^{-4F}\frac{\Omega}{\mathsf{c}}\Pi^2\Big)  \nonumber \\
&=\Big(
\Big(1+\frac{\Omega}{\mathsf{c}}A\Big)^2
-e^{-4F}\frac{\Omega^2}{\mathsf{c}^2}\Pi^2\Big)^{-1}
\Big(A\Big(1+\frac{\Omega}{\mathsf{c}}\Big)-e^{-4F}\frac{\Omega}{\mathsf{c}}\Pi^2\Big).
\end{align}
Hence the identity $\nabla_{\mu}T^{\mu 1}=\nabla_{\mu}T^{\mu 3}=0$ reduces
\begin{align}
&\partial_1P+(\epsilon+P)\partial_1G
-(\epsilon+P)U^0U_2\frac{\partial_1\Omega}{\mathsf{c}} 
=0,  \nonumber \\
&\partial_3P+(\epsilon+P)\partial_3G
-(\epsilon+P)U^0U_2\frac{\partial_3\Omega}{\mathsf{c}}=0. \label{2.7}
\end{align}
On the other hand we see that the other Euler equations
$$\nabla_{\mu}T^{\mu 0}=0,\qquad \nabla_{\mu}T^{\mu 2}=0 $$
hold automatically, for we are assuming $\partial_0\Omega=\partial_2\Omega=0$.\\

Therefore, defining the `relativistic enthalpy density' $u$ by
\begin{equation}
u=f^{u}(\rho):=\mathsf{c}^2\int_0^{\rho}\frac{dP}{\epsilon+P}=
\int_0^{\rho}\frac{dP}{\rho+P/\mathsf{c}^2}, \label{Def.u}
\end{equation}
we have
\begin{equation}
\frac{u}{\mathsf{c}^2}+G=\mbox{Const.},
\end{equation}
while $\rho >0$, provided that $\Omega$ is constant while $\rho>0$. \\

If we consider the dust: $P=0$ in the region $\rho>=0$, where $\Omega$ is supposed to be constant, then 
\eqref{2.7} reads $\partial_1G=\partial_3G=0$ so that we set
\begin{equation}
G=\mbox{Const.}
\end{equation}

\subsection{Einstein equations}

In order to write down the Einstein equations
\begin{equation}
R_{\mu\nu}-\frac{1}{2}g_{\mu\nu}R=
\frac{8\pi\mathsf{G}}{\mathsf{c}^4}T_{\mu\nu}
\end{equation}
or
\begin{equation}
R_{\mu\nu}=\frac{8\pi\mathsf{G}}{\mathsf{c}^4}
(T_{\mu\nu}-\frac{1}{2}g_{\mu\nu}T),
\end{equation}
where $R$ stands for $g^{\alpha\beta}R_{\alpha\beta}$ and
$T$ stands for $g_{\alpha\beta}T^{\alpha\beta}$, let us compute $U_{\mu}$,
$T_{\mu\nu}$ and $T$. The result is as following:

\begin{align*}
&U_0=e^{2F-G}\Big(1+\frac{\Omega}{\mathsf{c}}A\Big), \\
&U_1=U_3=0, \\
&U_2=e^{-G+2F}\Big(
A\Big(1+\frac{\Omega}{\mathsf{c}}A\Big)-e^{-4F}\frac{\Omega}{\mathsf{c}}\Pi^2\Big) ;
\end{align*}
\begin{align*}
T_{00}&=(\epsilon+P)e^{4F-2G}\Big(1+\frac{\Omega}{\mathsf{c}}A\Big)^2
-Pe^{2F}, \\
T_{02}&=T_{20}= \\
&=(\epsilon+P)e^{4F-2G}\Big(1+\frac{\Omega}{\mathsf{c}}A\Big)
\Big(A\Big(1+\frac{\Omega}{\mathsf{c}}A\Big)-
e^{-4F}\frac{\Omega}{\mathsf{c}}\Pi^2\Big)
-Pe^{2F}A, \\
T_{11}&=T_{33}=Pe^{-2F+2K}, \\
T_{22}&=
(\epsilon+P)e^{-2G+4F}
\Big(A\Big(1+\frac{\Omega}{\mathsf{c}}A\Big)-e^{-4F}\frac{\Omega}{\mathsf{c}}\Pi^2\Big)^2
-P(e^{2F}A^2-e^{-2F}\Pi^2)
\end{align*}
and other $T_{\mu\nu}$'s are zero; Let us note that 
$$
T=(\epsilon+P)e^{-2G}\Big[e^{2F}\Big(1+\frac{\Omega}{\mathsf{c}}A\Big)^2
-e^{-2F}\frac{\Omega^2}{\mathsf{c}^2}\Pi^2\Big]
-4P
$$
reduces to
\begin{equation}
T=\epsilon-3P. \label{TeP}
\end{equation}\\

\subsection{Lewis metric expression by $f, k, l, m$}

In order to calculate the Ricci tensor
$$ R_{\mu\nu}=\partial_{\alpha}\Gamma_{\mu\nu}^{\alpha}
-\partial_{\nu}\Gamma_{\mu\alpha}^{\alpha}
+\Gamma_{\mu\nu}^{\alpha}\Gamma_{\alpha\beta}^{\beta}
-\Gamma_{\mu\alpha}^{\beta}\Gamma_{\nu\beta}^{\alpha}
$$
 and to write down explicitly the Einstein equations, it is convenient to write the metric as
\begin{equation}
ds^2=f\mathsf{c}^2dt^2-2k\mathsf{c}dtd\phi-ld\phi^2
-e^m(d\varpi^2+dz^2), \label{25}
\end{equation}
that is, to put
\begin{equation}
f=e^{2F},\quad k=-e^{2F}A, \quad m=2(-F+K),\quad l=-e^{2F}A^2+
e^{-2F}\Pi^2.
\end{equation}

The expression \eqref{25} is called `{\bf Lewis metric}' after \cite{Lewis}\\

First let us note the identity
\begin{equation}
\Pi^2=fl+k^2. \label{27}
\end{equation}
and that the factor $e^{-G}=U^0$ is given by
\begin{equation}
(U^0)^{-2}=e^{2G}=f-2\frac{\Omega}{\mathsf{c}}k
-\frac{\Omega^2}{\mathsf{c}^2}l.
\end{equation}

The the components of the metric and the Christoffel symbols describing by $f, k, l, m$ are given in {\bf Appendix 2}.

The components of the Ricci tensor are as following (other $R_{\mu\nu}$ are zero ):
\begin{subequations}
\begin{align}
\frac{2e^m}{\Pi}R_{00}&=\partial_1\frac{\partial_1f}{\Pi}+\partial_3\frac{\partial_3f}{\Pi}+
\frac{1}{\Pi^3}f\Sigma, \label{30a}\\
-\frac{2e^m}{\Pi}R_{02}&=-\frac{2e^m}{\Pi}R_{20}=
\partial_1\frac{\partial_1k}{\Pi}+\partial_3\frac{\partial_3k}{\Pi}+\frac{1}{\Pi^3}k\Sigma, \label{30b} \\
-\frac{2e^m}{\Pi}R_{22}&=
\partial_1\frac{\partial_1l}{\Pi}+\partial_3\frac{\partial_3l}{\Pi}+\frac{1}{\Pi^3}l\Sigma, \label{30c}\\
2R_{11}&=-\partial_1^2m-\partial_3^2m-2\frac{\partial_1^2\Pi}{\Pi}
+\frac{1}{\Pi}[(\partial_1m)(\partial_1\Pi)-(\partial_3m)(\partial_3\Pi)] + \nonumber \\
&+\frac{1}{\Pi^2}[(\partial_1f)(\partial_1l)+(\partial_1k)^2], \label{30d}\\
2R_{33}&=-\partial_1^2m-\partial_3^2m-2\frac{\partial_3^2\Pi}{\Pi}
-\frac{1}{\Pi}[(\partial_1m)(\partial_1\Pi)-(\partial_3m)(\partial_3\Pi)] + \nonumber \\
&+\frac{1}{\Pi^2}[(\partial_3f)(\partial_3l)+(\partial_3k)^2] \label{30e}\\
2R_{13}&=2R_{31}=-2\frac{\partial_1\partial_3\Pi}{\Pi}+
\frac{1}{\Pi}[(\partial_3m)(\partial_1\Pi)+(\partial_1m)(\partial_3\Pi)]
\nonumber \\
&+\frac{1}{2\Pi^2}[(\partial_1f)(\partial_3l)+(\partial_1l)(\partial_3f)
+2(\partial_1k)(\partial_3k)].\label{30f}
\end{align}
\end{subequations}

Here we have introduced the quantity $\Sigma$ defined by
\begin{equation}
\Sigma:=(\partial_1f)(\partial_1l)+(\partial_3f)(\partial_3l)+
(\partial_1k)^2+(\partial_3k)^2,
\end{equation}
and during the calculations we have used the identities
\begin{equation}
\Gamma_{j\alpha}^{\alpha}=\partial_jm+\frac{\partial_j\Pi}{\Pi} \qquad
\mbox{for}\quad j=1,3.
\end{equation}
Note that 
$$\Sigma=\sum_{j=1,3}(\partial_jf)(\partial_jl)+(\partial_jk)^2 $$
with
\begin{equation}
(\partial_jf)(\partial_jl)+(\partial_jk)^2 = e^{4F}(\partial_jA)^2
-4(\partial_jF)^2\Pi^2 +4\Pi(\partial_j\Pi)(\partial_jF). \label{SigmaF}
\end{equation}

Now the components of the 4-velocity vector are:
\begin{align*}
&U^0=e^{-G}=\Big(f-
2\frac{\Omega}{\mathsf{c}}k-\frac{\Omega^2}{\mathsf{c}^2}l\Big)^{-1/2},
\quad U^1=U^3=0,\quad U^2=\frac{\Omega}{\mathsf{c}}U^0; \\
&U_0=\Big(f-\frac{\Omega}{\mathsf{c}}k\Big)U^0,\quad
U_1=U_3=0,\quad
U_2=-\Big(k+\frac{\Omega}{\mathsf{c}}l\Big)U^0.
\end{align*}

The Einstein equations are
\begin{equation}
R_{\mu\nu}=\frac{8\pi\mathsf{G}}{\mathsf{c}^4}S_{\mu\nu}, \quad\mbox{where}\quad
S_{\mu\nu}:=T_{\mu\nu}-\frac{1}{2}g_{\mu\nu}T.
\end{equation}

The components $S_{\mu\nu}$ turn out to be as following (other
$S_{\mu\nu}$ are zero) :
\begin{subequations}
\begin{align}
&S_{00}=
\frac{1}{2}(\epsilon+P)e^{-2G}
\Big[\Big(f-\frac{\Omega}{\mathsf{c}}k\Big)^2+\frac{\Omega^2}{\mathsf{c}^2}\Pi^2\Big]
+Pf, \label{Ta}\\
&S_{02}=S_{20}=
\frac{1}{2}(\epsilon+P)e^{-2G}
\Big[-kf-2\frac{\Omega}{\mathsf{c}}fl+\frac{\Omega^2}{\mathsf{c}^2}kl\Big]-Pk, \label{Tb}\\
&S_{22}=
\frac{1}{2}(\epsilon+P)e^{-2G}\Big[\Pi^2+
\Big(k+\frac{\Omega}{\mathsf{c}}l\Big)^2\Big]-Pl, \label{Tc}\\
&S_{11}=S_{33}=\frac{e^m}{2}(\epsilon -P). \label{Td}
\end{align}
\end{subequations}

Recall \eqref{TeP}.\\

Thus the full set of Einstein equations is:
\begin{subequations}
\begin{align}
&R_{00}=\frac{8\pi\mathsf{G}}{\mathsf{c}^4}S_{00}, \label{36a}\\
&R_{02}=\frac{8\pi\mathsf{G}}{\mathsf{c}^4}S_{02}, \label{36b}\\
&R_{22}=\frac{8\pi\mathsf{G}}{\mathsf{c}^4}S_{22}, \label{36c}\\
&R_{11}=\frac{8\pi\mathsf{G}}{\mathsf{c}^4}S_{11}, \label{36d} \\
&R_{33}=\frac{8\pi\mathsf{G}}{\mathsf{c}^4}S_{33}, \label{36e}\\
&R_{13}=0 \label{36f}
\end{align}
\end{subequations}

. 

\section{Equations for $F, A, \Pi, K$}

{\bf 1) } As for \eqref{30a} we see that 
\begin{align}
\frac{2e^m}{f}R_{00}&=
\frac{\Pi}{2f}\Big[\partial_1\Big(\frac{1}{\Pi}\partial_1f\Big)
+\partial_3\Big(\frac{1}{\Pi}\partial_3f\Big)+\frac{1}{\Pi^3}f\Sigma\Big]= \nonumber \\
&=\partial_1^2F+\partial_3^2F+\sum_{j=1,3}
\Big[\frac{1}{\Pi}(\partial_j\Pi)(\partial_jF)+\frac{e^{4F}}{2\Pi^2}(\partial_jA)^2\Big].
\end{align}
On the other hand, \eqref{36a} reads
\begin{equation}
\frac{2e^m}{f}R_{00}=
\frac{4\pi\mathsf{G}}{\mathsf{c}^4}\frac{e^m}{f}
\Big[(\epsilon+P)e^{-2G}((f-\frac{\Omega}{\mathsf{c}}k)^2+\frac{\Omega^2}{\mathsf{c}^2}\Pi^2)+2Pf\Big]
\end{equation}
so that
\begin{align}
&\frac{2e^m}{f}R_{00}= \nonumber \\
&=\partial_1^2F+\partial_3^2F+\sum_{j=1,3}
\Big[\frac{1}{\Pi}(\partial_j\Pi)(\partial_jF)+\frac{e^{4F}}{2\Pi^2}(\partial_jA)^2\Big]= \nonumber \\
&=\frac{4\pi\mathsf{G}}{\mathsf{c}^4}\frac{e^m}{f}
\Big[(\epsilon+P)e^{-2G}((f-\frac{\Omega}{\mathsf{c}}k)^2+\frac{\Omega^2}{\mathsf{c}^2}\Pi^2)+2Pf\Big]. \label{N.1}
\end{align}
But we see
$$e^{-2G}\Big(\Big(f-\frac{\Omega}{\mathsf{c}}k\Big)^2+\frac{\Omega^2}{\mathsf{c}^2}\Pi^2\Big)= f\frac{e^{2F}\Big(1+\frac{\Omega}{\mathsf{c}}A\Big)^2+e^{-2F}\frac{\Omega^2}{\mathsf{c}^2}\Pi^2}{e^{2F}\Big(1+\frac{\Omega}{\mathsf{c}}A\Big)^2-e^{-2F}\frac{\Omega^2}{\mathsf{c}^2}\Pi^2}.
$$
Thus \eqref{N.1} reads
\begin{align}
&\partial_1^2F+\partial_3^2F+\sum_{j=1,3}
\Big[\frac{1}{\Pi}(\partial_j\Pi)(\partial_jF)+\frac{e^{4F}}{2\Pi^2}(\partial_jA)^2\Big] = \nonumber \\
&=\frac{4\pi\mathsf{G}}{\mathsf{c}^4}
e^{2(-F+K)}
\Big[(\epsilon+P)
\frac{e^{2F}\Big(1+\frac{\Omega}{\mathsf{c}}A\Big)^2+e^{-2F}\frac{\Omega^2}{\mathsf{c}^2}\Pi^2}{e^{2F}\Big(1+\frac{\Omega}{\mathsf{c}}A\Big)^2-e^{-2F}\frac{\Omega^2}{\mathsf{c}^2}\Pi^2}
+2P\Big] \label{N.2}
\end{align}\\

{\bf 2)} As for \eqref{30b}  we see that
\begin{align}
\frac{2e^m}{f}R_{02}&=
-\frac{\Pi}{f}\Big[\partial_1\Big(\frac{1}{\Pi}\partial_1k\Big)+
\partial_3\Big(\frac{1}{\Pi}\partial_3k\Big)+\frac{1}{\Pi^3}k\Sigma \Big] = \nonumber \\
&=
\partial_1^2A+\partial_3^2A+\sum_{j=1,3}\Big[-\frac{1}{\Pi}(\partial_j\Pi)(\partial_jA)+
4(\partial_jF)(\partial_jA) + \nonumber \\
&+\Big(\frac{e^{4F}}{\Pi^2}(\partial_jA)^2
+2\partial_j^2F+2\frac{1}{\Pi}(\partial_j\Pi)(\partial_jF)\Big)\cdot A\Big] . \label{N.3}
\end{align}
On the other hand, \eqref{36b} reads
\begin{equation}
\frac{2e^m}{f}R_{02}=\frac{8\pi\mathsf{G}}{\mathsf{c}^4}\frac{e^m}{f}\Big[
(\epsilon+P)e^{-2G}\Big(-kf-2\frac{\Omega}{\mathsf{c}}fl+\frac{\Omega2}{\mathsf{c}^2}kl\Big)
-2kP\Big]. 
\end{equation}
But we see
$$e^{-2G}\Big(-kf-2\frac{\Omega}{\mathsf{c}}fl+\frac{\Omega^2}{\mathsf{c}^2}kl\Big)=
f\Big[
A-2\frac{\Omega}{\mathsf{c}}\frac{e^{-2F}\Pi^2}{e^{2F}\Big(1+\frac{\Omega}{\mathsf{c}}A\Big)^2-e^{-2F}\frac{\Omega^2}{\mathsf{c}^2}\Pi^2}
\Big].
$$
Thus \eqref{N.3} reads
\begin{align}
&\partial_1^2A+\partial_3^2A+\sum_{j=1,3}\Big[-\frac{1}{\Pi}(\partial_j\Pi)(\partial_jA)+
4(\partial_jF)(\partial_jA) + \nonumber \\
&+\Big(\frac{e^{4F}}{\Pi^2}(\partial_jA)^2
+2\partial_j^2F+2\frac{1}{\Pi}(\partial_j\Pi)(\partial_jF)\Big)\cdot A\Big] = \nonumber \\
&=\frac{8\pi\mathsf{G}}{\mathsf{c}^4}
e^{2(-F+K)}\Big[(\epsilon+P)
\Big(
A-2\frac{\Omega}{\mathsf{c}}\frac{e^{-2F}\Pi^2}{e^{2F}\Big(1+\frac{\Omega}{\mathsf{c}}A\Big)^2-e^{-2F}\frac{\Omega^2}{\mathsf{c}^2}\Pi^2}
\Big)+2AP\Big]. \label{N.4}
\end{align}

Taking \eqref{N.4}$- 2A\times$\eqref{N.2}, we get
\begin{align}
&\frac{2e^m}{f}(R_{02}-2A\cdot R_{00})= \nonumber \\
&=\partial_1^2A+\partial_3^2A+\sum_{j=1,3}\Big[-\frac{1}{\Pi}(\partial_j\Pi)(\partial_jA)+
4(\partial_jF)(\partial_jA) \Big]= \nonumber \\
&=-\frac{16\pi\mathsf{G}}{\mathsf{c}^4}
e^{2(-F+K)}(\epsilon+P)
\frac{e^{-2F}\frac{\Omega}{\mathsf{c}}\Pi^2\Big(1+\frac{\Omega}{\mathsf{c}}A\Big)}{e^{2F}\Big(1+\frac{\Omega}{\mathsf{c}}A\Big)^2-e^{-2F}\frac{\Omega^2}{\mathsf{c}^2}\Pi^2}. \label
{N.5}
\end{align}\\

{\bf 3)} We have the identity
\begin{equation}
lS_{00}-2kS_{02}
-fS_{22}=2P\Pi^2, \label{N.6}
\end{equation}
which can be verified from \eqref{Ta}\eqref{Tb}\eqref{Tc} thanks to \eqref{27}. 
On the other hand,
we have the identity
\begin{equation}
\frac{e^m}{\Pi}(lR_{00}-2kR_{02}-fR_{22})=
\partial_1^2\Pi+\partial_3^2\Pi,\label{N.7}
\end{equation}
which can be verified from \eqref{30a}\eqref{30b}\eqref{30c} thanks to \eqref{27}.
Hence \eqref{N.6} and \eqref{N.7} leads us to the equation
\begin{equation}
\partial_1^2\Pi+\partial_3^2\Pi=\frac{16\pi\mathsf{G}}{\mathsf{c}^4}
e^{2(-F+K)}P\Pi, \label{N.8}
\end{equation}
if \eqref{36a}\eqref{36b}\eqref{36c} hold.

\begin{Remark} If we consider the dust for which $P=0$, then \eqref{N.8} says that $\Pi(\varpi, z)$ is a harmonic function and we can assume $\Pi=\varpi$ by conformal change of coordinates.
See \cite[p. 26]{Islam}. (This was first used by 
\cite{Weyl} and generalized to the present case by
\cite{Lewis}.)  But it is not the case when $P\not=0$. 
\end{Remark}

Summing up, we have the following

\begin{Proposition}
The set of equations \eqref{36a} \eqref{36b} \eqref{36c}
implies that
\begin{subequations}
\begin{align}
&\frac{\partial^2F}{\partial\varpi^2}+\frac{\partial^2F}{\partial z^2}
+\frac{1}{\Pi}\Big(\frac{\partial F}{\partial\varpi}\frac{\partial\Pi}{\partial\varpi}
+\frac{\partial F}{\partial z}\frac{\partial \Pi}{\partial z}\Big)+
\frac{e^{4F}}{2\Pi^2}\Big[\Big(\frac{\partial A}{\partial\varpi}\Big)^2+
\Big(\frac{\partial A}{\partial z}\Big)^2\Big] = \nonumber \\
&=\frac{4\pi\mathsf{G}}{\mathsf{c}^4}
e^{2(-F+K)}
\Big[(\epsilon+P)
\frac{e^{2F}\Big(1+\frac{\Omega}{\mathsf{c}}A\Big)^2+e^{-2F}\frac{\Omega^2}{\mathsf{c}^2}\Pi^2}{e^{2F}\Big(1+\frac{\Omega}{\mathsf{c}}A\Big)^2-e^{-2F}\frac{\Omega^2}{\mathsf{c}^2}\Pi^2}
+2P\Big], \label{N.9a} \\
&\frac{\partial^2A}{\partial\varpi^2}+\frac{\partial^2A}{\partial z^2}
-\frac{1}{\Pi}\Big(\frac{\partial\Pi}{\partial\varpi}\frac{\partial A}{\partial\varpi}+\frac{\partial\Pi}{\partial z}\frac{\partial A}{\partial z}\Big)+
4\Big(\frac{\partial F}{\partial\varpi}\frac{\partial A}{\partial\varpi}+
\frac{\partial F}{\partial z}\frac{\partial A}{\partial z}\Big) = \nonumber \\
&=-\frac{16\pi\mathsf{G}}{\mathsf{c}^4}
e^{2(-F+K)}(\epsilon+P)
\frac{e^{-2F}\frac{\Omega}{\mathsf{c}}\Pi^2\Big(1+\frac{\Omega}{\mathsf{c}}A\Big)}{e^{2F}\Big(1+\frac{\Omega}{\mathsf{c}}A\Big)^2-e^{-2F}\frac{\Omega^2}{\mathsf{c}^2}\Pi^2}, \label{N.9b} \\
&\frac{\partial^2\Pi}{\partial \varpi^2}
+\frac{\partial^2\Pi}{\partial z^2}=\frac{16\pi\mathsf{G}}{\mathsf{c}^4}
e^{2(-F+K)}P\Pi. \label{N.9c}
\end{align}
\end{subequations}
\end{Proposition}

Inversely we can show that \eqref{N.9a}\eqref{N.9b}\eqref{N.9c} imply 
\eqref{36a}\eqref{36b}\eqref{36c}.
In fact, if we put
$$Q_{\mu\nu}:=R_{\mu\nu}-\frac{8\pi\mathsf{G}}{\mathsf{c}^4}S_{\mu\nu},
$$
then \eqref{36a}\eqref{36b}\eqref{36c} claim that $Q_{00}=Q_{02}=Q_{22}=0$.
On the other hand \eqref{N.9a}\eqref{N.9b}\eqref{N.9c} are nothing but
\begin{align*}
&Q_{00}=0, \\
&-2AQ_{00}
+Q_{02} =0, \\
&lQ_{00}-2kQ_{02}-fQ_{22}=0.
\end{align*}
Here we see 
\begin{equation}
\Delta:=\det
\begin{bmatrix}
1 & 0 & 0 \\
-2A & 1 & 0 \\
l & -2k & -f 
\end{bmatrix} =-f=-e^{2F}\not= 0.
\end{equation}

Thus we can claim

\begin{Proposition}\label{PropositionNB}
The set of equations
\eqref{N.9a}, \eqref{N.9b}, \eqref{N.9c} implies 
\eqref{36a}, \eqref{36b}, \eqref{36c}.
\end{Proposition}

Proof. $\Delta \not=0$ guarantees
 $Q_{00}=Q_{02}=Q_{22}=0$.
This means \eqref{36a},\eqref{36b},\eqref{36c}. $\square$. \\

{\bf 4)} It follows from \eqref{30d}\eqref{30e} that
\begin{align*}
2(R_{11}-R_{33})&=-\frac{2}{\Pi}(\partial_1^2\Pi-\partial_3^2\Pi)+
\frac{2}{\Pi}[(\partial_1m)(\partial_1\Pi)-(\partial_3m)(\partial_3\Pi)]+ \\
&+\frac{1}{\Pi^2}
[(\partial_1f)(\partial_1l)+(\partial_1k)^2
-(\partial_3f)(\partial_3l)-(\partial_3k)^2)].
\end{align*}
Therefore \eqref{Td} implies
\begin{align}
-\frac{2}{\Pi}(\partial_1^2\Pi-\partial_3^2\Pi)+&
\frac{2}{\Pi}[(\partial_1m)(\partial_1\Pi)-(\partial_3m)(\partial_3\Pi)]+ \nonumber \\
&+\frac{1}{\Pi^2}
[(\partial_1f)(\partial_1l)+(\partial_1k)^2
-(\partial_3f)(\partial_3l)-(\partial_3k)^2]=0.\label{N.53}
\end{align}
Recall \eqref{SigmaF} and $m=-2F+2K$. So,
 \eqref{N.53} reads
\begin{align}
(\partial_1\Pi)(\partial_1K)-(\partial_3\Pi)(\partial_3K)&=
\frac{1}{2}(\partial_1^2\Pi-\partial_3^2\Pi)+
\Pi[(\partial_1F)^2-(\partial_3F)^2] + \nonumber \\
&-\frac{e^{4F}}{4\Pi}
[(\partial_1A)^2-(\partial_3A)^2].
\end{align}\\

{\bf 5)} \eqref{36f} and \eqref{30f} imply
that
\begin{align*}
-2\frac{\partial_1\partial_3\Pi}{\Pi}&+
\frac{1}{\Pi}[(\partial_3m)(\partial_1\Pi)+
(\partial_1m)(\partial_3\Pi)] + \\
&+\frac{1}{2\Pi^2}
[(\partial_1f)(\partial_3l)+(\partial_1l)(\partial_3f)+2(\partial_1k)(\partial_3k)]=0.
\end{align*}
But we see that
\begin{align}
&(\partial_1f)(\partial_3l)+(\partial_1l)(\partial_3f)+2(\partial_1k)(\partial_3k)= \nonumber \\
&=-8(\partial_1F)(\partial_3F)\Pi^2+
4\Pi[(\partial_1F)(\partial_3\Pi)+(\partial_3F)(\partial_1\Pi)]+
2e^{4F}(\partial_1A)(\partial_3A ).
\end{align}

Therefore, since $m=-2F+2K$, \eqref{36f} reads
\begin{equation}
(\partial_3\Pi)(\partial_1K)+(\partial_1\Pi)(\partial_3K)=
\partial_1\partial_3\Pi+2\Pi(\partial_1F)(\partial_3F)-
\frac{e^{4F}}{2\Pi}(\partial_1A)(\partial_3A ).
\end{equation}\\

Thus we can claim

\begin{Proposition}
The set of equations \eqref{36d} \eqref{36e} \eqref{36f}, provided \eqref{Td}, implies
\begin{subequations}
\begin{align}
(\partial_1\Pi)(\partial_1K)-(\partial_3\Pi)(\partial_3K)&=
\frac{1}{2}(\partial_1^2\Pi-\partial_3^2\Pi)+
\Pi[(\partial_1F)^2-(\partial_3F)^2] + \nonumber \\
&-\frac{e^{4F}}{4\Pi}
[(\partial_1A)^2-(\partial_3A)^2], \label{Ya} \\
(\partial_3\Pi)(\partial_1K)+(\partial_1\Pi)(\partial_3K)=&
\partial_1\partial_3\Pi+2\Pi(\partial_1F)(\partial_3F)-
\frac{e^{4F}}{2\Pi}(\partial_1A)(\partial_3A ). \label{Yb}
\end{align}
\end{subequations}
\end{Proposition}

However the inverse is doubtful, namely, we are not sure that the system of equations \eqref{Ya} and \eqref{Yb} can recover both \eqref{36d} and  \eqref{36e}, separately, when $\Omega$ is not a constant. \\

Be that as it may, we suppose the following assumption:\\

{\bf (B2):} {\it It holds that }
\begin{equation}
\Big(\frac{\partial\Pi}{\partial\varpi}\Big)^2+\Big(\frac{\partial \Pi}{\partial z}\Big)^2
\not=0. 
\end{equation}\\

Under this assumption the set of equations \eqref{Ya}\eqref{Yb} is equivalent to 
\begin{subequations}
\begin{align}
&\partial_1K=((\partial_1\Pi)^2+(\partial_3\Pi)^2)^{-1}\Big[
(\partial_1\Pi)\mathrm{RH}\eqref{Ya}+(\partial_3\Pi)\mathrm{RH}\eqref{Yb}\Big], \label{Yc} \\
&\partial_3K=((\partial_1\Pi)^2+(\partial_3\Pi)^2)^{-1}\Big[
-(\partial_3\Pi)\mathrm{RH}\eqref{Ya}+(\partial_1\Pi)\mathrm{RH}\eqref{Yb}\Big], \label{Yd}
\end{align}
\end{subequations}
where RH\eqref{Ya}, RH\eqref{Yb}  stand for the right-hand sides of 
\eqref{Ya}, \eqref{Yb}, respectively.

We can claim
\begin{Proposition}\label{Prop4}
Suppose the assumption {\bf (B2)} and that either $\Omega$ is a constant
on the considered domain or $\rho=P=0$, namely vacuum,  on the considered domain.
 Then the set of equations
\eqref{Ya},\eqref{Yb} implies \eqref{36d},\eqref{36e},\eqref{36f}, provided that the equation
\eqref{Td} and the set of equations \eqref{N.9a},\eqref{N.9b}, \eqref{N.9c} hold. 
\end{Proposition}

Proof.  Since the set of equations \eqref{Ya}\eqref{Yb} is equivalent to the set of equations
\eqref{36d}$-$\eqref{36e}, \eqref{36f}, we have to prove \eqref{36d}$+$\eqref{36e}, namely, we have to prove
\begin{equation}
\frac{1}{2}(R_{11}+R_{33})=\frac{8\pi\mathsf{G}}{\mathsf{c}^4}S_{11}=
\frac{8\pi\mathsf{G}}{\mathsf{c}^4}S_{33}=\frac{4\pi\mathsf{G}}{\mathsf{c}^4}e^m(\epsilon-P).
\label{EqProp4}
\end{equation}

So, we consider 
$$-\frac{1}{2}(R_{11}+R_{33})\cdot (\Pi_1^2+\Pi_3^2)=\clubsuit, $$
where
$$\clubsuit:=(\Pi_1^2+\Pi_3^2)\Big[-\triangle F+\triangle K+
\frac{\triangle\Pi}{2\Pi}-\frac{\Sigma}{4\Pi^2}\Big].$$
Here $\Pi_j, (j=1,3), \triangle\Pi, \triangle K, \triangle F$ stand for
$\partial_j\Pi$,
$$\frac{\partial^2 \Pi}{\partial\varpi^2}+\frac{\partial^2\Pi}{\partial z^2}, \quad
\frac{\partial^2 K}{\partial\varpi^2}+\frac{\partial^2K}{\partial z^2}, \quad
\frac{\partial^2F}{\partial\varpi^2}+\frac{\partial^2F}{\partial z^2}. $$

Differentiating \eqref{Yc} by $\varpi$ and \eqref{Yd} by $z$, we have
\begin{align*}
&(\Pi_1^2+\Pi_3^2)\triangle K= \\
&=-(\Pi_1K_1+\Pi_3K_3)\triangle\Pi+
\frac{1}{2}(\Pi_1\partial_1\triangle\Pi+\Pi_3\partial_3\triangle\Pi)+ \\
&=(\Pi_1^2-\Pi_3^2)(F_1^2-F_3^2)+4\Pi_1\Pi_3F_1F_3+
2\Pi(\Pi_1F_1+\Pi_3F_3)\triangle F + \\
&+\frac{e^{4F}}{4\Pi^2}\Big[(\Pi_1^2-\Pi_3^2+4\Pi(-\Pi_1F_1+\Pi_3F_3))(A_1^2-A_3^2) + \\
&+(2\Pi_1\Pi_3-4\Pi(\Pi_1F_3+\Pi_3F_1))(2A_1A_3) + \\
&-2\Pi(\Pi_1A_1+\Pi_3A_3)\triangle A \Big]
\end{align*}
after tedious calculations. Here $ F_j, A_j, ( j=1,3,)$ stand for
$\partial_jF, \partial_jA$, and
$\triangle A$ means
$ \displaystyle
\frac{\partial^2A}{\partial\varpi^2}+\frac{\partial^2A}{\partial z^2} $.

Denoting by $[S\mathrm{a}], [S\mathrm{b}]$ the right-hand sides of the equations
\eqref{N.9a}, \eqref{N.9b}, we eliminate $\triangle F, \triangle A$. The result is
\begin{align*}
&(\Pi_1^2+\Pi_3^2)\triangle K= \\
&=-(\Pi_1K_1+\Pi_3K_3)\triangle\Pi+
\frac{1}{2}(\Pi_1\partial_1\triangle\Pi+\Pi_3\partial_3\triangle\Pi)+ \\
&-(\Pi_1^2+\Pi_3^2)(F_1^2+F_3^2)-
\frac{e^{4F}}{4\Pi^2}(\Pi_1^2+\Pi_3^2)(A_1^2+A_3^2) + \\
&+2\Pi (\Pi_1F_1+\Pi_3F_3)[S\mathrm{a}]
-\frac{e^{4F}}{2\Pi}(\Pi_1A_1+\Pi_3A_3)[S\mathrm{b}].
\end{align*}

Let us consider the case of constant $\Omega$. Put
\begin{align}
f':=e^{2G}&=e^{2F}\Big(1+\frac{\Omega}{\mathsf{c}}A\Big)^2
-e^{-2F}\frac{\Omega^2}{\mathsf{c}^2}\Pi^2 \nonumber \\
&=f-2\frac{\Omega}{\mathsf{c}}k-\frac{\Omega^2}{\mathsf{c}}l.
\end{align}
Then, using the identity
$$\partial_j\triangle\Pi=
\Big(2(-F_j+K_j)+\frac{\partial_jP}{P}+\frac{\Pi_j}{\Pi}\Big)\triangle\Pi $$
with
$$
\partial_jP=(\epsilon+P)\frac{1}{2f'}\Big(-\partial_jf+
2\frac{\Omega}{\mathsf{c}}\partial_jk+\frac{\Omega^2}{\mathsf{c}^2}\partial_jl\Big),
$$
 provided that $\Omega$ is a constant, and using
\begin{align*}
[S\mathrm{a}]&=\frac{4\pi\mathsf{G}}{\mathsf{c}^4}\frac{e^m}{f}\Big[(\epsilon+P)\frac{1}{f'}((f-\frac{\Omega}{\mathsf{c}}k)^2+\frac{\Omega^2}{\mathsf{c}^2}\Pi^2))+2Pf\Big], \\
[S\mathrm{b}]&=\frac{16\pi\mathsf{G}}{\mathsf{c}^4}\frac{e^m\Pi^2}{f^2f'}
(\epsilon+P)\frac{\Omega}{\mathsf{c}}\Big(\frac{\Omega}{\mathsf{c}}k-f\Big),
\end{align*}
we can deduce
$$\clubsuit =(\Pi_1^2+\Pi_3^2)(-\epsilon+P)\frac{4\pi\mathsf{G}}{\mathsf{c}^4}e^m.$$
Here we have used the identity
$$\Sigma=e^{4F}(A_1^2+A_3^2)-4\Pi^2(F_1^2+F_3^2)+4\Pi(\Pi_1F_1+\Pi_3F_3).$$
This gives the desired equation \eqref{EqProp4}.\\

When the vacuum is considered, we have $\triangle\Pi=[S\mathrm{a}]=[S\mathrm{b}]=0$, and we see $\clubsuit =0$. This completes the proof. $\square$.\\

{\bf 6)} Summing up, we get the system of equations
\begin{subequations}
\begin{align}
&\frac{\partial^2F}{\partial\varpi^2}+\frac{\partial^2F}{\partial z^2}
+\frac{1}{\Pi}\Big(\frac{\partial F}{\partial\varpi}\frac{\partial\Pi}{\partial\varpi}
+\frac{\partial F}{\partial z}\frac{\partial \Pi}{\partial z}\Big)+
\frac{e^{4F}}{2\Pi^2}\Big[\Big(\frac{\partial A}{\partial\varpi}\Big)^2+
\Big(\frac{\partial A}{\partial z}\Big)^2\Big] = \nonumber \\
&=\frac{4\pi\mathsf{G}}{\mathsf{c}^4}
e^{2(-F+K)}
\Big[(\epsilon+P)
\frac{e^{2F}\Big(1+\frac{\Omega}{\mathsf{c}}A\Big)^2+e^{-2F}\frac{\Omega^2}{\mathsf{c}^2}\Pi^2}{e^{2F}\Big(1+\frac{\Omega}{\mathsf{c}}A\Big)^2-e^{-2F}\frac{\Omega^2}{\mathsf{c}^2}\Pi^2}
+2P\Big], \label{N.EQa} \\
&\frac{\partial^2A}{\partial\varpi^2}+\frac{\partial^2A}{\partial z^2}
-\frac{1}{\Pi}\Big(\frac{\partial\Pi}{\partial\varpi}\frac{\partial A}{\partial\varpi}+\frac{\partial\Pi}{\partial z}\frac{\partial A}{\partial z}\Big)+
4\Big(\frac{\partial F}{\partial\varpi}\frac{\partial A}{\partial\varpi}+
\frac{\partial F}{\partial z}\frac{\partial A}{\partial z}\Big) = \nonumber \\
&=-\frac{16\pi\mathsf{G}}{\mathsf{c}^4}
e^{2(-F+K)}(\epsilon+P)
\frac{e^{-2F}\frac{\Omega}{\mathsf{c}}\Pi^2\Big(1+\frac{\Omega}{\mathsf{c}}A\Big)}{e^{2F}\Big(1+\frac{\Omega}{\mathsf{c}}A\Big)^2-e^{-2F}\frac{\Omega^2}{\mathsf{c}^2}\Pi^2}, \label{N.EQb} \\
&\frac{\partial^2\Pi}{\partial \varpi^2}
+\frac{\partial^2\Pi}{\partial z^2}=\frac{16\pi\mathsf{G}}{\mathsf{c}^4}
e^{2(-F+K)}P\Pi, \label{N.EQc} \\
&\frac{\partial\Pi}{\partial \varpi}
\frac{\partial K}{\partial\varpi}-
\frac{\partial\Pi}{\partial z}\frac{\partial K}{\partial z}=
\frac{1}{2}\Big(\frac{\partial^2\Pi}{\partial\varpi^2}
-\frac{\partial^2\Pi}{\partial z^2}\Big)+
\Pi\Big[\Big(\frac{\partial F}{\partial\varpi}\Big)^2
-\Big(\frac{\partial F}{\partial z}\Big)^2\Big] + \nonumber \\
&-\frac{e^{4F}}{4\Pi}
\Big[\Big(\frac{\partial A}{\partial\varpi}\Big)^2
-\Big(\frac{\partial A}{\partial z}\Big)^2\Big], \label{N.EQd} \\
&\frac{\partial\Pi}{\partial z}
\frac{\partial K}{\partial\varpi}
+\frac{\partial\Pi}{\partial \varpi}\frac{\partial K}{\partial z}=
\frac{\partial^2\Pi}{\partial \varpi\partial z}
+2\Pi\frac{\partial F}{\partial\varpi}
\frac{\partial F}{\partial z}-
\frac{e^{4F'}}{2\Pi}
\frac{\partial A}{\partial\varpi}\frac{\partial A}{\partial z} , \label{N.EQe} \\
&G=
\Big( =\frac{1}{2}\log\Big[ e^{2F}\Big(1+\frac{\Omega}{\mathsf{c}}A\Big)^2-
e^{-2F}\frac{\Omega^2}{\mathsf{c}^2}\Pi^2\Big]\quad \Big) = \nonumber \\
&=-\frac{u}{\mathsf{c}^2}+\mbox{Const.}.\label{N.EQf}
\end{align}
\end{subequations}

Here $u, P, \epsilon=\mathsf{c}^2\rho$ are given functions of $\rho$, and $G$ is
determined by $F, A, \Pi, \Omega$ through \eqref{18}.\\

\begin{Theorem}
Suppose the assumption {\bf (B2)} and that either $\Omega$ is a constant on the considered domain or $\rho=P=0$ on the considered domain.  Then the system of equations
\eqref{N.EQa} $\sim$ \eqref{N.EQe} is equivalent to the system of Einstein equations
\eqref{36a} $\sim$ \eqref{36f}.
\end{Theorem}

{\bf 7)}  Now we have a question of the consistency of the first order system of equations for $K$. In order that there exists $K$ which satisfies \eqref{N.EQd}, \eqref{N.EQe}, or there exists $K$ which satisfies 
\begin{align}
\frac{\partial K}{\partial\varpi}&=\Big(
\Big(\frac{\partial\Pi}{\partial\varpi}\Big)^2+\Big(\frac{\partial \Pi}{\partial z}\Big)^2\Big)^{-1}
\Big[\frac{\partial\Pi}{\partial\varpi}\mbox{RH\eqref{N.EQd}}
+\frac{\partial\Pi}{\partial z}\mbox{RH\eqref{N.EQe}}\Big], \label{Z.Ka} \\
\frac{\partial K}{\partial z}&=\Big(
\Big(\frac{\partial\Pi}{\partial\varpi}\Big)^2+\Big(\frac{\partial \Pi}{\partial z}\Big)^2\Big)^{-1}
\Big[-\frac{\partial\Pi}{\partial z}\mbox{RH\eqref{N.EQd}}
+\frac{\partial\Pi}{\partial\varpi}\mbox{RH\eqref{N.EQe}}\Big], \label{Z.Kb}
\end{align}
where RH\eqref{N.EQd}, RH\eqref{N.EQe} stand for the right-hand sides of 
\eqref{N.EQd}, \eqref{N.EQe}, provided the assumption {\bf (B2)}, it is necessary that the `{\bf consistency condition}' 
\begin{equation}
\frac{\partial \tilde{K}_1}{\partial z}=\frac{\partial \tilde{K}_3}{\partial\varpi} \label{Z.Cons}
\end{equation}
holds, where $\tilde{K}_1, \tilde{K}_3$ stand for the right-hand sides of 
\eqref{Z.Ka}, \eqref{Z.Kb}. 

It is claimed in \cite[Section 4.2, p.56]{Islam} that it is the case when $P=0$ and $\Omega$ is a constant in the considered domain, if we take $\Pi=\varpi$. Actually we can claim
the following: 

\begin{Proposition}\label{PropN5}
Suppose the assumption {\bf (B2)} and 
let $\mathfrak{D}$ be a domain which is the union of a domain $\mathfrak{D}_0$ on which $\Omega$ is a constant and a domain $\mathfrak{D}_1$ on which $\rho=P=0$, namely, a  vacuum domain. Let $K$ be arbitrarily given, and let $F, A, \Pi$ satisfy \eqref{N.EQa}, \eqref{N.EQb}, \eqref{N.EQc}, and \eqref{N.EQf} for $\rho=0, \Omega$ and this given $K$. Denote by 
$\tilde{K}_1, \tilde{K}_3$ the right-hand sides of \eqref{Z.Ka}, \eqref{Z.Kb}, respectively evaluated by these $F, A, \Pi$. Then it holds on $\mathfrak{D}$ that 
\begin{align}
\frac{\partial\tilde{K}_1}{\partial z}-
\frac{\partial\tilde{K}_3}{\partial\varpi}&=
\frac{16\pi \mathsf{G}}{\mathsf{c}^4}
e^{2(-F+K)}P\Pi\Big[\Big(\frac{\partial\Pi}{\partial\varpi}\Big)^2+
\Big(\frac{\partial\Pi}{\partial z}\Big)^2\Big]^{-1}\times \nonumber \\
&\times \Big[\Big(\frac{\partial K}{\partial\varpi}-\tilde{K}_1\Big)
\frac{\partial\Pi}{\partial z}-
\Big(\frac{\partial K}{\partial z}-\tilde{K}_3\Big)
\frac{\partial\Pi}{\partial\varpi}\Big].\label{N3.39}
\end{align}
\end{Proposition}

Proof. By a tedious calculation, we get
\begin{align*}
\frac{\partial\tilde{K}_1}{\partial z}-\frac{\partial\tilde{K}_3}{\partial\varpi}&=
-(\Pi_3\tilde{K}_1-\Pi_1\tilde{K}_3)
(\Pi_1^2+\Pi_3^2)^{-1}\cdot \triangle\Pi + \\
&+(\Pi_1^2+\Pi_3^2)^{-1}Z
\end{align*}
with
\begin{align*}
Z:=&\frac{1}{2}(-\Pi_1\partial_3\triangle\Pi+\Pi_3\partial_1\triangle\Pi)+ \\
&+2\Pi(-\Pi_1F_3+\Pi_3F_1)[S\mathrm{a}]
+\frac{e^{4F}}{2\Pi}(\Pi_1A_3-\Pi_3A_1)[S\mathrm{b}],
\end{align*}
where $\Pi_j, F_j, A_j, (j=1,3,) \triangle\Pi$ stand for $\partial_j\Pi, \partial_jF, \partial_jA, \partial_1^2\Pi+\partial_3^2\Pi$ respectively, and $[S\mathrm{a}], [S\mathrm{b}] $ stand for the right-hand sides of the equations \eqref{N.EQa}, \eqref{N.EQb}. 

Consider $Z$ on the domain $\mathfrak{D}_0$ on which $\Omega$ is a constant. We can write
\begin{align*}
Z=&\frac{1}{2}(-\Pi_1\partial_3\triangle\Pi+\Pi_3\partial_1\triangle\Pi) + \\
&+(-\Pi_1f_3+\Pi_3f_1)\Big(\frac{\Pi}{f}[S\mathrm{a}]-\frac{k}{2\Pi}[S\mathrm{b}]\Big) + \\
&+(-\Pi_1k_3+\Pi_3k_1)\frac{f}{2\Pi}[S\mathrm{b}],
\end{align*}
where $f_j, k_j$ stand for $\partial_jf, \partial_jk$. We have
\begin{align*} 
&\frac{\Pi}{f}[S\mathrm{a}]-\frac{k}{2\Pi}[S\mathrm{b}]=
\frac{8\pi\mathsf{G}}{\mathsf{c}^4}e^m\Big[
(\epsilon+P)\frac{\Pi}{2f^2f'}\Big(f^2+\frac{\Omega^2}{\mathsf{c}^2}(\Pi^2-k^2)\Big)
+\frac{\Pi P}{f}\Big], \\
&[S\mathrm{b}]=\frac{16\pi\mathsf{G}}{\mathsf{c}^4}\frac{e^m\Pi^2}{f^2f'}
(\epsilon+P)\frac{\Omega}{\mathsf{c}}\Big(\frac{\Omega}{\mathsf{c}}k-f\Big),
\end{align*}
where $f'=e^{2G}=f-2\frac{\Omega}{\mathsf{c}}k-\frac{\Omega^2}{\mathsf{c}^2}l$.  Using
the identity
$$\partial_j\triangle\Pi=
\Big(2(-F_j+K_j)+\frac{\partial_jP}{P}+\frac{\Pi_j}{\Pi}\Big)\triangle\Pi
$$
with
$$
\partial_jP=(\epsilon+P)\frac{1}{2f'}\Big(-f_j+
2\frac{\Omega}{\mathsf{c}}k_j+\frac{\Omega^2}{\mathsf{c}^2}\partial_jl\Big) ,
$$
provided that $\Omega$ is a constant, we can deduce that
$$Z=(-\Pi_1\partial_3K+\Pi_3\partial_1K)\triangle\Pi.$$
Here we have used the identity
$$\partial_jl=\frac{1}{f}\Big(2\Pi\Pi_j-\frac{\Pi^2-k^2}{f}f_j-2kk_j\Big),$$
which can be derived from the identity
$\Pi^2=fl+k^2$.
This implies the desired identity. 

On the domain on which $\rho=P=0$, we have $\triangle\Pi=[S\mathrm{a}]=
[S\mathrm{b}]=0$ so that $Z=0$ and
$$\frac{\partial\tilde{K}_1}{\partial z}-\frac{\partial \tilde{K}_3}{\partial \varpi}=0.$$
This completes the proof. $\square$\\

Therefore, as a conclusion of Proposition \ref{PropN5}, if $K$ satisfies \eqref{Z.Ka} \eqref{Z.Kb},
then the consistency condition 
\begin{equation}
\frac{\partial}{\partial z}\mbox{RH\eqref{Z.Ka}}=\frac{\partial}{\partial\varpi}\mbox{RH\eqref{Z.Kb}}
\end{equation}
holds, since 
$$\frac{\partial K}{\partial\varpi}=\tilde{K}_1,\qquad
\frac{\partial K}{\partial z}=\tilde{K}_3.
$$ Of course this conclusion in itself is a vicious circular argument of no use. However the following argument is useful: 

\begin{Theorem} \label{ThN2}
Let us consider a bounded domain $\mathfrak{D}$ which is a union of a domain 
$\mathfrak{D}_0$ on which $\Omega$ is a constant and a domain
$\mathfrak{D}_1$ on which $\rho=P=0$.
Suppose that $K \in C^1(\bar{\mathfrak{D}})$ is given and that $F, A, \Pi, \rho \in C^3(\bar{\mathfrak{D}})$ satisfy \eqref{N.EQa},\eqref{N.EQb},\eqref{N.EQc} and \eqref{N.EQf} with $\rho, \Omega$ and this $K$. 
Suppose that the assumption  {\bf (B2)} holds on $\bar{\mathfrak{D}}$. Let us denote by
$\tilde{K}_1,\tilde{K}_3$ the right-hand sides of
\eqref{Z.Ka},\eqref{Z.Kb}, respectively, evaluated by these
$F, A, \Pi$. (They are $C^1$-functions on $\bar{\mathfrak{D}}$.)
Put
\begin{equation}
\tilde{K}(\varpi, z):=K_O+
\int_0^z\tilde{K}_3(0,z')dz'+
\int_0^{\varpi}\tilde{K}_1(\varpi',z)d\varpi'\label{N3.40}
\end{equation}
for $(\varpi, z) \in \mathfrak{D}$. Here $K_O$ is a constant. If $\tilde{K}=K$, then $K$
satisfies 
\begin{equation}
\frac{\partial K}{\partial\varpi}=\tilde{K}_1,\qquad
\frac{\partial K}{\partial z}=\tilde{K}_3,\label{N3.41}
\end{equation}
that is, the equations \eqref{N.EQd}\eqref{N.EQe} are satisfied.
\end{Theorem}

 Proof.  Suppose that $\tilde{K}=K$. It follows from \eqref{N3.40} with $\tilde{K}=K$ that
\begin{align}
\frac{\partial K}{\partial \varpi}(\varpi, z)&=\tilde{K}_1(\varpi, z), \label{N3.42} \\
\frac{\partial K}{\partial z}(\varpi, z)&=
\tilde{K}_3(0,z)+\int_0^{\varpi}
\frac{\partial\tilde{K}_1}{\partial z}(\varpi', z)d\varpi'.\label{N3.43}
\end{align}
Put
\begin{equation}
L(\varpi, z):=\frac{\partial \tilde{K}_1}{\partial z}-
\frac{\partial\tilde{K}_3}{\partial\varpi},
\end{equation}
which is a continuous function on $\bar{\mathfrak{D}}$. Then \eqref{N3.43} reads
\begin{equation}
\frac{\partial K}{\partial z}(\varpi, z)=
\tilde{K}_3(\varpi, z)+
\int_0^{\varpi}
L(\varpi', z)d\varpi'.\label{N3.45}
\end{equation}
Now therefore \eqref{N3.39} of Proposition \ref{PropN5} reads
\begin{equation}
L(\varpi, z)=-
\frac{16\pi\mathsf{G}}{\mathsf{c}^4}
e^{2(-F+K)}P\Pi\Big[\Big(\frac{\partial\Pi}{\partial\varpi}\Big)^2+
\Big(\frac{\partial\Pi}{\partial z}\Big)^2\Big]^{-1}\frac{\partial\Pi}{\partial\varpi}
\int_0^{\varpi}L(\varpi',z)d\varpi'.
\end{equation}
Since the function
$$
\frac{16\pi\mathsf{G}}{\mathsf{c}^4}
e^{2(-F+K)}P\Pi\Big[\Big(\frac{\partial\Pi}{\partial\varpi}\Big)^2+
\Big(\frac{\partial\Pi}{\partial z}\Big)^2\Big]^{-1}\frac{\partial\Pi}{\partial\varpi}
$$
is bounded on the compact $\bar{\mathfrak{D}}$, the Gronwall's argument implies that
$L(\varpi, z)=0$ on $\mathfrak{D}$ so that \eqref{N3.45} reads
\begin{equation}
\frac{\partial K}{\partial z}(\varpi, z)=\tilde{K}_3(\varpi, z).\label{N3.47}
\end{equation}
Thus \eqref{N3.42} and \eqref{N3.47} complete the proof. $\square$\\

\section{Equations for $F', A', \Pi, K'$}

In this section we shall use the variables $f', k', l', m'$ defined by 
\begin{equation}
f':=f-2\frac{\Omega}{\mathsf{c}}k-\frac{\Omega^2}{\mathsf{c}^2}l,\quad
k':=k+\frac{\Omega}{\mathsf{c}}l,\quad l'=l,\quad m'=m
\end{equation}
$$\Leftrightarrow$$
\begin{equation}
f=f'+2\frac{\Omega}{\mathsf{c}}k'-\frac{\Omega^2}{\mathsf{c}^2}l',
\quad k=k'-\frac{\Omega}{\mathsf{c}}l',\quad l=l',\quad m=m'.
\end{equation}
and define $F', A', K'$ by
\begin{equation}
e^{2F'}=f'(=e^{2G}),\quad k'=-e^{2F'}A',\quad m'=2(-F'+K'). \label{3.21}
\end{equation}
Let us note the identities
\begin{equation}
\Pi^2=f'l'+(k')^2, \label{3.23}
\end{equation}
and
\begin{equation}
l'=-e^{2F'}(A')^2+e^{-2F'}\Pi^2. \label{3.24}
\end{equation}\\

{\bf 1)}
Now it follows from \eqref{Ta}\eqref{Tb}\eqref{Tc} that
\begin{equation}
S_{00}+2\frac{\Omega}{\mathsf{c}}S_{02}
+\frac{\Omega^2}{\mathsf{c}^2}S_{22}=
\frac{1}{2}(\epsilon+3P)(f-2\frac{\Omega}{\mathsf{c}}k-\frac{\Omega^2}{\mathsf{c}^2}l). \label{3.27}
\end{equation}

In order to calculate the corresponding
$$R_{00}+2\frac{\Omega}{\mathsf{c}}R_{02}
+\frac{\Omega^2}{\mathsf{c}^2}R_{22},$$
we use the variables  $f', k', l', m'$.

Then it follows from \eqref{30a}\eqref{30b}\eqref{30c} that
\begin{equation}
\frac{2e^m}{\Pi}\Big(
R_{00}+2\frac{\Omega}{\mathsf{c}}R_{02}+\frac{\Omega^2}{\mathsf{c}^2}R_{22}\Big)=
\partial_1\frac{\partial_1f'}{\Pi}+\partial_3\frac{\partial_3f'}{\Pi}+\frac{f'}{\Pi^3}\Sigma
+[D1], \label{3.28}
\end{equation}
where 
\begin{equation}
[D1]=[D1,1]+[D1,3]
\end{equation}
with
\begin{equation}
[D1,j]:=\frac{4}{\Pi}\frac{\partial_j\Omega}{\mathsf{c}}\partial_jk'+
\partial_j\Big(\frac{2}{\Pi}\frac{\partial_j\Omega}{\mathsf{c}}\Big)k'-
\frac{2}{\Pi}\Big(\frac{\partial_j\Omega}{\mathsf{c}}\Big)^2l'.
\end{equation}

Thus, \eqref{3.27}, \eqref{3.28} give
\begin{equation}
\partial_1\frac{\partial_1f'}{\Pi}+\partial_3\frac{\partial_3f'}{\Pi}+\frac{f'}{\Pi^3}\Sigma
+[D1]=\frac{8\pi\mathsf{G}}{\mathsf{c}^4}e^m\frac{f'}{\Pi}(\epsilon +3P), \label{3.30}
\end{equation}
provided \eqref{36a}\eqref{36b}\eqref{36c}.\\

On the other hand, using \eqref{3.23} and \eqref{3.21}, we can verify that
\begin{equation}
(\partial_jf')(\partial_jl')+(\partial_jk')^2=
2\Pi\partial_j\Pi\frac{\partial_jf'}{f'}-
\Big(\frac{\Pi\partial_jf'}{f'}\Big)^2
+(f'\partial_jA')^2.
\end{equation}
So, putting
\begin{equation}
\Sigma':=(\partial_1 f')(\partial_1l')+(\partial_3f')(\partial_3l')+
(\partial_1k')^2+(\partial_3k')^2,
\end{equation}
we have
\begin{equation}
\Sigma'=\sum_{j=1,3}
2\Pi\partial_j\Pi\frac{\partial_jf'}{f'}-
\Big(\frac{\Pi\partial_jf'}{f'}\Big)^2
+(f'\partial_jA')^2.
\end{equation}
But we have
\begin{equation}
(\partial_jf)(\partial_jl)+(\partial_jk)^2=
(\partial_jf')(\partial_jl')+(\partial_jk')^2+[W1,j]
\end{equation}
with
\begin{equation}
[W1, j]:=2\frac{\partial_j\Omega}{\mathsf{c}}(k'\partial_jl'-l'\partial_jk')+
\Big(\frac{\partial_j\Omega}{\mathsf{c}}\Big)^2(l')^2.
\end{equation}

Therefore
we have
\begin{equation}
\Sigma=\Sigma'+[W1],
\end{equation}
where
\begin{equation}
[W1]=[W1,1]+[W1,3].
\end{equation}

Therefore, \eqref{3.30}, which can be written as
$$
\sum_{j=1,3}\frac{\Pi}{f'}\partial_j\frac{\partial_jf'}{\Pi}+
\frac{\Sigma}{\Pi^2}+\frac{\Pi}{f'}[D1]=\frac{8\pi\mathsf{G}}{\mathsf{c}^4}e^m(\epsilon +3P),
$$
reads
\begin{align}
&\frac{2e^m}{f'}
\Big(
R_{00}+2\frac{\Omega}{\mathsf{c}}R_{02}+\frac{\Omega^2}{\mathsf{c}^2}R_{22}\Big) = \nonumber \\
&=
\partial_1^2F'+\partial_3^2F'+\sum_{j=1,3}\Big[
\frac{(\partial_j\Pi)(\partial_jF')}{\Pi}+
\frac{e^{4F'}}{2\Pi^2}(\partial_jA')^2\Big]
+\Pi e^{-2F'}[D1]+\frac{[W1]}{2\Pi^2}
= \nonumber \\
&= \frac{4\pi\mathsf{G}}{\mathsf{c}^4}e^{2(-F'+K')}(\epsilon+3P). \label{3.35}
\end{align}

When $\Omega$ is a constant, then $[D1]=[W1]=0$, and the equation \eqref{3.35}
 is nothing but \cite[(1.34a)]{Meinel}. \\

{\bf 2)}
Let us note the identity
\begin{equation}
\Big(k+\frac{\Omega}{\mathsf{c}}l\Big)f-\Big(f+\frac{\Omega^2}{\mathsf{c}^2}l\Big)k
-\frac{\Omega}{\mathsf{c}}\Big(f-\frac{\Omega}{\mathsf{c}}\Big)l=0.
\end{equation}

So, we can verify that
\begin{equation}
\Big(k+\frac{\Omega}{\mathsf{c}}l\Big)S_{00}+
\Big(f+\frac{\Omega^2}{\mathsf{c}^2}l\Big)S_{02}
+\frac{\Omega}{\mathsf{c}}\Big(f-\frac{\Omega}{\mathsf{c}}k\Big)S_{22}=0.
\end{equation}

Therefore

\begin{equation}
-\frac{2e^m}{\Pi}
\Big[\Big(k+\frac{\Omega}{\mathsf{c}}l\Big)R_{00}+
\Big(f+\frac{\Omega^2}{\mathsf{c}^2}l\Big)R_{02}+
\frac{\Omega}{\mathsf{c}}\Big(f-\frac{\Omega}{\mathsf{c}}k\Big)R_{22}\Big] =  0 \label{3.18}
\end{equation}

Let us calculate the right-hand side of \eqref{3.18}, say RH\eqref{3.18}.

We get
\begin{equation}
\mbox{RH\eqref{3.18}} = \partial_1\frac{1}{\Pi}(f'\partial_1k'-k'\partial_1f')+
\partial_3\frac{1}{\Pi}(f'\partial_3k'-k'\partial_3f') -[D2],
\end{equation}
where
\begin{equation}
[D2]=[D2,1]+[D2,3]
\end{equation}
with
\begin{equation}
[D2,j]:=\frac{1}{\Pi}\frac{\partial_j^2\Omega}{\mathsf{c}}\Big(\Pi^2+\frac{2\Omega}{\mathsf{c}}k'l'\Big)
-\frac{1}{\Pi}\frac{\partial_j\Omega}{\mathsf{c}}l'\partial_jf'.
\end{equation}

Since
\begin{equation}
f'\partial_jk'-k'\partial_jf'=-e^{4F'}\partial_jA',
\end{equation}
 \eqref{3.18} turns out to be
\begin{equation}
-\partial_1\frac{1}{\Pi}e^{4F'}\partial_1A'
-\partial_3\frac{1}{\Pi}e^{4F'}\partial_3A' = [D2], \label{3.26}
\end{equation}
or
\begin{align}
&\frac{2e^m}{(f')^2}
\Big[\Big(k+\frac{\Omega}{\mathsf{c}}l\Big)R_{00}+
\Big(f+\frac{\Omega^2}{\mathsf{c}^2}l\Big)R_{02}+
\frac{\Omega}{\mathsf{c}}\Big(f-\frac{\Omega}{\mathsf{c}}k\Big)R_{22}\Big] = \nonumber \\
&=\partial_1^2A'+\partial_3^2A'+
\sum_{j=1,3}\Big(-\frac{1}{\Pi}(\partial_j\Pi)(\partial_jA')+4(\partial_jF')(\partial_jA')\Big)
+[D2]
\nonumber \\
&=0.
\end{align}\\

When $\Omega$ is a constant, then $[D2]=0$ and \eqref{3.26} is nothing but \cite[(1.34b)]{Meinel}.\\

{\bf 3)}
We have the identity
\begin{equation}
lS_{00}-2kS_{02}
-fS_{22}=2P\Pi^2, \label{37}
\end{equation}
which can be verified from \eqref{Ta}\eqref{Tb}\eqref{Tc} thanks to \eqref{27}. 
On the other hand,
we have the identity
\begin{equation}
\frac{e^m}{\Pi}(lR_{00}-2kR_{02}-fR_{22})=
\partial_1^2\Pi+\partial_3^2\Pi,\label{32}
\end{equation}
which can be verified from \eqref{30a}\eqref{30b}\eqref{30c} thanks to \eqref{27}.
Hence \eqref{32} and \eqref{37} leads us to the equation
\begin{equation}
\partial_1^2\Pi+\partial_3^2\Pi=\frac{16\pi\mathsf{G}}{\mathsf{c}^4}
e^mP\Pi, \label{38}
\end{equation}
if \eqref{36a}\eqref{36b}\eqref{36c} hold.

Summing up, we have the following

\begin{Proposition}
The set of equations \eqref{36a} \eqref{36b} \eqref{36c}
implies that
\begin{subequations}
\begin{align}
&\frac{\partial^2F'}{\partial\varpi^2}+\frac{\partial^2F'}{\partial z^2}
+\frac{1}{\Pi}\Big(\frac{\partial F'}{\partial\varpi}\frac{\partial\Pi}{\partial\varpi}
+\frac{\partial F'}{\partial z}\frac{\partial \Pi}{\partial z}\Big)+
\frac{e^{4F'}}{2\Pi^2}\Big[\Big(\frac{\partial A'}{\partial\varpi}\Big)^2+
\Big(\frac{\partial A'}{\partial z}\Big)^2\Big] + \nonumber \\
&+\Pi e^{-2F'}[D1]+\frac{[W1]}{2\Pi^2} = \frac{4\pi\mathsf{G}}{\mathsf{c}^4}e^{2(-F'+K')}(\epsilon+3P), \label{EAa} \\
&\frac{\partial^2A'}{\partial\varpi^2}+\frac{\partial^2A'}{\partial z^2}
-\frac{1}{\Pi}\Big(\frac{\partial\Pi}{\partial\varpi}\frac{\partial A'}{\partial\varpi}+\frac{\partial\Pi}{\partial z}\frac{\partial A'}{\partial z}\Big)+
4\Big(\frac{\partial F'}{\partial\varpi}\frac{\partial A'}{\partial\varpi}+
\frac{\partial F'}{\partial z}\frac{\partial A'}{\partial z}\Big) + \nonumber \\
&+[D2]=0, \label{EAb} \\
&\frac{\partial^2\Pi}{\partial \varpi^2}
+\frac{\partial^2\Pi}{\partial z^2}=\frac{16\pi\mathsf{G}}{\mathsf{c}^4}
e^{2(-F'+K')}P\Pi. \label{EAc}
\end{align}
\end{subequations}
\end{Proposition}

Inversely we can show that \eqref{EAa}\eqref{EAb}\eqref{EAc} imply 
\eqref{36a}\eqref{36b}\eqref{36c}.
In fact, if we put
$$Q_{\mu\nu}:=R_{\mu\nu}-\frac{8\pi\mathsf{G}}{\mathsf{c}^4}S_{\mu\nu},
$$
then \eqref{36a}\eqref{36b}\eqref{36c} claim that $Q_{00}=Q_{02}=Q_{22}=0$.
On the other hand \eqref{EAa}\eqref{EAb}\eqref{EAc} are nothing but
\begin{align*}
&Q_{00}+2\frac{\Omega}{\mathsf{c}}Q_{02}+\frac{\Omega^2}{\mathsf{c}^2}Q_{22}=0, \\
&\Big(k+\frac{\Omega}{\mathsf{c}}l\Big)Q_{00}
+\Big(f+\frac{\Omega^2}{\mathsf{c}^2}l\Big)Q_{02}+
\frac{\Omega}{\mathsf{c}}\Big(f-\frac{\Omega}{\mathsf{c}}k\Big)Q_{22}=0, \\
&lQ_{00}-2kQ_{02}-fQ_{22}=0.
\end{align*}
Here we see 
\begin{align}
\Delta&:=\det
\begin{bmatrix}
l & -2k & -f \\
k+\frac{\Omega}{\mathsf{c}}l & f+\frac{\Omega^2}{\mathsf{c}^2}l & \frac{\Omega}{\mathsf{c}}\Big(f-\frac{\Omega}{\mathsf{c}}k\Big) \\
1 & 2\frac{\Omega}{\mathsf{c}} & \frac{\Omega^2}{\mathsf{c}^2} 
\end{bmatrix} \nonumber \\
&=\Big(2k+\frac{\Omega}{\mathsf{c}}l\Big)^2\frac{\Omega^2}{\mathsf{c}^2}
+f(f-\frac{4\Omega}{\mathsf{c}}k-\frac{2\Omega^2}{\mathsf{c}^2}l\Big) \nonumber \\
&=(f')^2=e^{4F'}\not= 0.
\end{align}

Thus we can claim

\begin{Proposition}\label{PropositionB}
The set of equations
\eqref{EAa}, \eqref{EAb}, \eqref{EAc} implies 
\eqref{36a}, \eqref{36b}, \eqref{36c}.
\end{Proposition}

Proof. $\Delta \not=0$ guarantees
 $Q_{00}=Q_{02}=Q_{22}=0$.
This means \eqref{36a},\eqref{36b},\eqref{36c}. $\square$. \\

{\bf 4)} It follows from \eqref{30d}\eqref{30e} that
\begin{align*}
2(R_{11}-R_{33})&=-\frac{2}{\Pi}(\partial_1^2\Pi-\partial_3^2\Pi)+
\frac{2}{\Pi}[(\partial_1m)(\partial_1\Pi)-(\partial_3m)(\partial_3\Pi)]+ \\
&+\frac{1}{\Pi^2}
[(\partial_1f)(\partial_1l)+(\partial_1k)^2
-(\partial_3f)(\partial_3l)-(\partial_3k)^2)].
\end{align*}
Therefore \eqref{Td} implies
\begin{align}
-\frac{2}{\Pi}(\partial_1^2\Pi-\partial_3^2\Pi)+&
\frac{2}{\Pi}(\partial_1m\partial_1\Pi-\partial_3m\partial_3\Pi)+ \nonumber \\
&+\frac{1}{\Pi^2}
((\partial_1f)(\partial_1l)+(\partial_1k)^2
-(\partial_3f)(\partial_3l)-(\partial_3k)^2)=0.\label{53}
\end{align}
Recall 
that we have
\begin{align}
(\partial_jf)(\partial_jl)+(\partial_jk)^2&=
(\partial_jf')(\partial_jl')+(\partial_jk')^2 +[W1, j] \nonumber\\
&=-\Big(\frac{\partial_jf'}{f'}\Big)^2\Pi^2+
2\frac{\partial_jf'}{f'}\Pi\partial_j\Pi+
(f')^2(\partial_jA')^2 +[W1, j]\nonumber \\
&=-4(\partial_jF')^2\Pi^2
+4(\partial_jF')\Pi(\partial_j\Pi)+
e^{4F'}(\partial_jA')^2+[W1, j].
\end{align}

Therefore, since $m=m'=-2F'+2K'$, \eqref{53} reads
\begin{align}
(\partial_1\Pi)(\partial_1K')-(\partial_3\Pi)(\partial_3K')&=
\frac{1}{2}(\partial_1^2\Pi-\partial_3^2\Pi)+
\Pi[(\partial_1F')^2-(\partial_3F')^2] + \nonumber \\
&-\frac{e^{4F'}}{4\Pi}
[(\partial_1A')^2-(\partial_3A')^2]-\frac{1}{4\Pi}([W1, 1]-[W1, 3]),
\end{align}
which is nothing but \cite[(1.35a)]{Meinel} when $\Omega$ is a constant so that
 $[W1, j]=0$.\\

{\bf 5)} \eqref{36f} and \eqref{30f} imply
that
\begin{align*}
-2\frac{\partial_1\partial_3\Pi}{\Pi}&+
\frac{1}{\Pi}((\partial_3m)(\partial_1\Pi)+
(\partial_1m)(\partial_3\Pi)) + \\
&+\frac{1}{2\Pi^2}
((\partial_1f)(\partial_3l)+(\partial_1l)(\partial_3f)+2(\partial_1k)(\partial_3k))=0.
\end{align*}
But we see that
\begin{align}
&(\partial_1f)(\partial_3l)+(\partial_1l)(\partial_3f)+2(\partial_1k)(\partial_3k)= \nonumber \\
&=(\partial_1f')(\partial_3l')+(\partial_1l')(\partial_3f')+2(\partial_1k')(\partial_3k') 
+2[W2] \nonumber \\
&=-8(\partial_1F')(\partial_3F')\Pi^2+
4\Pi[(\partial_1F')(\partial_3\Pi)+(\partial_3F')(\partial_1\Pi)]+
2e^{4F'}(\partial_1A')(\partial_3A' )+2[W2],
\end{align}
where
\begin{equation}
[W2]:=\frac{\partial_1\Omega}{\mathsf{c}}(k'\partial_3l'-l'\partial_3k')
+\frac{\partial_3\Omega}{\mathsf{c}}(k'\partial_1l'-l'\partial_1k')
+\frac{\partial_1\Omega}{\mathsf{c}}\frac{\partial_3\Omega}{\mathsf{c}}(l')^2.
\end{equation}

Therefore, since $m=m'=-2F'+2K'$, \eqref{36f} reads
\begin{align}
(\partial_3\Pi)(\partial_1K')+(\partial_1\Pi)(\partial_3K')=&
\partial_1\partial_3\Pi+2\Pi(\partial_1F')(\partial_3F')-
\frac{e^{4F'}}{2\Pi}(\partial_1A')(\partial_3A' )+ \nonumber \\
&-\frac{[W2]}{2\Pi},
\end{align}
which is nothing but \cite[(1.35b)]{Meinel}
when $\Omega$ is a constant so that $[W2]=0$.\\

Thus we can claim
\begin{Proposition}
The set of equations \eqref{36d} \eqref{36e} \eqref{36f}, provided \eqref{Td}, implies
\begin{subequations}
\begin{align}
(\partial_1\Pi)(\partial_1K')-(\partial_3\Pi)(\partial_3K')&=
\frac{1}{2}(\partial_1^2\Pi-\partial_3^2\Pi)+
\Pi[(\partial_1F')^2-(\partial_3F')^2] + \nonumber \\
&-\frac{e^{4F'}}{4\Pi}
[(\partial_1A')^2-(\partial_3A')^2]-\frac{1}{4\Pi}([W1, 1]-[W1, 3]), \label{Xa} \\
(\partial_3\Pi)(\partial_1K')+(\partial_1\Pi)(\partial_3K')=&
\partial_1\partial_3\Pi+2\Pi(\partial_1F')(\partial_3F')-
\frac{e^{4F'}}{2\Pi}(\partial_1A')(\partial_3A' )+ \nonumber \\
&-\frac{[W2]}{2\Pi}. \label{Xb}
\end{align}
\end{subequations}
\end{Proposition}

However the inverse is not so trivial. 

We are supposing the assumption {\bf (B2)}:
$$
\Big(\frac{\partial\Pi}{\partial\varpi}\Big)^2+\Big(\frac{\partial \Pi}{\partial z}\Big)^2
\not=0. 
$$
Under this assumption the set of equations \eqref{Xa}\eqref{Xb} is equivalent to 
\begin{subequations}
\begin{align}
&\partial_1K'=((\partial_1\Pi)^2+(\partial_3\Pi)^2)^{-1}\Big[
(\partial_1\Pi)\mathrm{RH}\eqref{Xa}+(\partial_3\Pi)\mathrm{RH}\eqref{Xb}\Big], \label{Xc} \\
&\partial_3K'=((\partial_1\Pi)^2+(\partial_3\Pi)^2)^{-1}\Big[
-(\partial_3\Pi)\mathrm{RH}\eqref{Xa}+(\partial_1\Pi)\mathrm{RH}\eqref{Xb}\Big], \label{Xd}
\end{align}
\end{subequations}
where RH\eqref{Xa}, RH\eqref{Xb}  stand for the right-hand sides of 
\eqref{Xa}, \eqref{Xb}, respectively.

We can claim
\begin{Proposition}\label{Prop6}
Suppose the assumption {\bf (B2)} and that $\Omega$ is a constant. Then the set of equations
\eqref{Xa},\eqref{Xb} implies \eqref{36d},\eqref{36e},\eqref{36f}, provided that the equation
\eqref{Td} and the set of equations \eqref{EAa},\eqref{EAb}, \eqref{EAc} hold (of course with $[D1]=[W1]=[D2]=0$). 
\end{Proposition}

Proof. Since the set of equations \eqref{Xa}\eqref{Xb} is equivalent to the set of equations
\eqref{36d}$-$\eqref{36e}, \eqref{36f}, we have to prove \eqref{36d}$+$\eqref{36e}, namely,
\begin{equation}
\frac{1}{2}(R_{11}+R_{33})=\frac{8\pi\mathsf{G}}{\mathsf{c}^4}S_{11}=
\frac{8\pi\mathsf{G}}{\mathsf{c}^4}S_{33}=\frac{4\pi\mathsf{G}}{\mathsf{c}^4}e^m(\epsilon-P).
\label{EqProp6}
\end{equation}

Differentiating \eqref{Xc} by $\varpi$ and \eqref{Xd} by $z$, we have
\begin{align*}
&(\Pi_1^2+\Pi_3^2)\triangle K'= \\
&=\Big[-(\Pi_1F_1'+\Pi_3F_3')+\frac{1}{2P}(\Pi_1\partial_1P+\Pi_3\partial_3P)+
\frac{1}{2\Pi}(\Pi_1^2+\Pi_3^2)\Big]\triangle\Pi + \\
&+(\Pi_1^2-\Pi_3^2)((F_1')^2-(F_3')^2)+4\Pi_1\Pi_3F_1'F_3'+
2\Pi(\Pi_1F_1'+\Pi_3F_3')\triangle F' + \\
&+\frac{e^{4F'}}{4\Pi^2}\Big[(\Pi_1^2-\Pi_3^2+4\Pi(-\Pi_1F_1'+\Pi_3F_3'))((A_1')^2-(A_3')^2) + \\
&+(2\Pi_1\Pi_3-4\Pi(\Pi_1F_3'+\Pi_3F_1'))(2A_1'A_3') + \\
&-2\Pi(\Pi_1A_1'+\Pi_3A_3')\triangle A' \Big]
\end{align*}
after tedious calculations. Here $\Pi_j, F_j', A_j', j=1,3,$ stand for
$\partial_j\Pi, \partial_jF', \partial_jA'$, and
$\triangle K', \triangle \Pi, \triangle F', \triangle A'$ mean
$$\frac{\partial^2 K'}{\partial\varpi^2}+\frac{\partial^2K'}{\partial z^2}, \quad
\frac{\partial^2\Pi}{\partial\varpi^2}+\frac{\partial^2\Pi}{\partial z^2}, \quad
\frac{\partial^2F'}{\partial\varpi^2}+\frac{\partial^2F'}{\partial z^2}, \quad
\frac{\partial^2A'}{\partial\varpi^2}+\frac{\partial^2A'}{\partial z^2}, $$
respectively. We have already used the equation \eqref{EAc} on $\triangle \Pi$.

Keeping in mind that
$$\partial_jP=-(\epsilon+P)\partial_jF',$$
we eliminate $\triangle\Pi, \triangle F', \triangle A'$ using the equations
\eqref{EAa}, \eqref{EAb}, \eqref{EAc}. Then, as the result of tedious calculations, we can get
\begin{align*}
&-\frac{1}{2}(R_{11}+R_{33})\cdot(\Pi_1^2+\Pi_3^2)=
(\Pi_1^2+\Pi_3^2)\Big[-\triangle F' +\triangle K' +\frac{\triangle\Pi}{2\Pi}-
\frac{\Sigma'}{4\Pi^2}\Big]= \\
&=(\Pi_1^2+\Pi_3^2)(-\epsilon +P)\cdot \frac{4\pi\mathsf{G}}{\mathsf{c}^4}e^{2m},
\end{align*}
that is, the desired equation \eqref{EqProp6}. 
Here we have used the identity
$$\Sigma'=4\Pi(\Pi_1F_1'+\Pi_3F_3')-4\Pi^2((F_1')^2+(F_3')^2)+e^{4F'}((A_1')^2+(A_3')^2).$$
This completes the proof. $\square$\\

{\bf 6)} Summing up, we get the system of equations

\begin{subequations}
\begin{align}
&\frac{\partial^2F'}{\partial\varpi^2}+\frac{\partial^2F'}{\partial z^2}
+\frac{1}{\Pi}\Big(\frac{\partial F'}{\partial\varpi}\frac{\partial\Pi}{\partial\varpi}
+\frac{\partial F'}{\partial z}\frac{\partial \Pi}{\partial z}\Big)+
\frac{e^{4F'}}{2\Pi^2}\Big[\Big(\frac{\partial A'}{\partial\varpi}\Big)^2+
\Big(\frac{\partial A'}{\partial z}\Big)^2\Big] + \nonumber \\
&+\Pi e^{-2F'}[D1]+\frac{[W1]}{2\Pi^2} =  \frac{4\pi\mathsf{G}}{\mathsf{c}^4}e^{2(-F'+K')}(\epsilon+3P), \label{EQa} \\
&\frac{\partial^2A'}{\partial\varpi^2}+\frac{\partial^2A'}{\partial z^2}
-\frac{1}{\Pi}\Big(\frac{\partial\Pi}{\partial\varpi}\frac{\partial A'}{\partial\varpi}+\frac{\partial\Pi}{\partial z}\frac{\partial A'}{\partial z}\Big)+
4\Big(\frac{\partial F'}{\partial\varpi}\frac{\partial A'}{\partial\varpi}+
\frac{\partial F'}{\partial z}\frac{\partial A'}{\partial z}\Big) +\nonumber \\
&+[D2]=0, \label{EQb} \\
&\frac{\partial^2\Pi}{\partial \varpi^2}
+\frac{\partial^2\Pi}{\partial z^2}=\frac{16\pi\mathsf{G}}{\mathsf{c}^4}
e^{2(-F'+K')}P\Pi, \label{EQc} \\
&\frac{\partial\Pi}{\partial \varpi}
\frac{\partial K'}{\partial\varpi}-
\frac{\partial\Pi}{\partial z}\frac{\partial K'}{\partial z}=
\frac{1}{2}\Big(\frac{\partial^2\Pi}{\partial\varpi^2}
-\frac{\partial^2\Pi}{\partial z^2}\Big)+
\Pi\Big[\Big(\frac{\partial F'}{\partial\varpi}\Big)^2
-\Big(\frac{\partial F'}{\partial z}\Big)^2\Big] + \nonumber \\
&-\frac{e^{4F'}}{4\Pi}
\Big[\Big(\frac{\partial A'}{\partial\varpi}\Big)^2
-\Big(\frac{\partial A'}{\partial z}\Big)^2\Big]-\frac{1}{4\Pi}([W1, 1]-[W1, 3]), \label{EQd} \\
&\frac{\partial\Pi}{\partial z}
\frac{\partial K'}{\partial\varpi}
+\frac{\partial\Pi}{\partial \varpi}\frac{\partial K'}{\partial z}=
\frac{\partial^2\Pi}{\partial \varpi\partial z}
+2\Pi\frac{\partial F'}{\partial\varpi}
\frac{\partial F'}{\partial z}-
\frac{e^{4F'}}{2\Pi}
\frac{\partial A'}{\partial\varpi}\frac{\partial A'}{\partial z}
 -\frac{[W2]}{2\Pi}, \label{EQe} \\
&F'=-\frac{u}{\mathsf{c}^2}+\mbox{Const.}.\label{EQf}
\end{align}
\end{subequations}

Here $u, P, \epsilon=\mathsf{c}^2\rho$ are given functions of $\rho$.\\

\begin{Theorem}
Suppose the assumption {\bf (B2)} and that $\Omega$ is a constant. Then the system of equations
\eqref{EQa} $\sim$ \eqref{EQe} is equivalent to the system of Einstein equations
\eqref{36a} $\sim$ \eqref{36f}.
\end{Theorem}

Hereafter we assume that $\Omega$ is a constant. Therefore the set of equations \eqref{EQa} - \eqref{EQf} is reduced to

\begin{subequations}
\begin{align}
&\frac{\partial^2F'}{\partial\varpi}+\frac{\partial^2F'}{\partial z}
+\frac{1}{\Pi}\Big(\frac{\partial F'}{\partial\varpi}\frac{\partial\Pi}{\partial\varpi}
+\frac{\partial F'}{\partial z}\frac{\partial \Pi}{\partial z}\Big)+
\frac{e^{4F'}}{2\Pi^2}\Big[\Big(\frac{\partial A'}{\partial\varpi}\Big)^2+
\Big(\frac{\partial A'}{\partial z}\Big)^2\Big] = \nonumber \\
&= \frac{4\pi\mathsf{G}}{\mathsf{c}^4}e^{2(-F'+K')}(\epsilon+3P)
, \label{CEQa} \\
&\frac{\partial^2A'}{\partial\varpi^2}+\frac{\partial^2A'}{\partial z^2}
-\frac{1}{\Pi}\Big(\frac{\partial\Pi}{\partial\varpi}\frac{\partial A'}{\partial\varpi}+\frac{\partial\Pi}{\partial z}\frac{\partial A'}{\partial z}\Big)+
4\Big(\frac{\partial F'}{\partial\varpi}\frac{\partial A'}{\partial\varpi}+
\frac{\partial F'}{\partial z}\frac{\partial A'}{\partial z}\Big) =0, \label{CEQb} \\
&\frac{\partial^2\Pi}{\partial \varpi^2}
+\frac{\partial^2\Pi}{\partial z^2}=\frac{16\pi\mathsf{G}}{\mathsf{c}^4}
e^{2(-F'+K')}P\Pi. \label{CEQc} \\
&\frac{\partial\Pi}{\partial \varpi}
\frac{\partial K'}{\partial\varpi}-
\frac{\partial\Pi}{\partial z}\frac{\partial K'}{\partial z}=
\frac{1}{2}\Big(\frac{\partial^2\Pi}{\partial\varpi^2}
-\frac{\partial^2\Pi}{\partial z^2}\Big)+
\Pi\Big[\Big(\frac{\partial F'}{\partial\varpi}\Big)^2
-\Big(\frac{\partial F'}{\partial z}\Big)^2\Big] + \nonumber \\
&-\frac{e^{4F'}}{4\Pi}
\Big[\Big(\frac{\partial A'}{\partial\varpi}\Big)^2
-\Big(\frac{\partial A'}{\partial z}\Big)^2\Big], \label{CEQd} \\
&\frac{\partial\Pi}{\partial z}
\frac{\partial K'}{\partial\varpi}
+\frac{\partial\Pi}{\partial \varpi}\frac{\partial K'}{\partial z}=
\frac{\partial^2\Pi}{\partial \varpi\partial z}
+2\Pi\frac{\partial F'}{\partial\varpi}
\frac{\partial F'}{\partial z}-
\frac{e^{4F'}}{2\Pi}
\frac{\partial A'}{\partial\varpi}\frac{\partial A'}{\partial z}
, \label{CEQe} \\
&F'=-\frac{u}{\mathsf{c}^2}+\mbox{Const.}.\label{CEQf}
\end{align}
\end{subequations}

As for the consistency condition of the first order system \eqref{CEQd} \eqref{CEQe} for
$K'$, or the system
\begin{subequations}
\begin{align}
\frac{\partial K'}{\partial\varpi}&=\Big[\Big(
\frac{\partial\Pi}{\partial\varpi}\Big)^2+
\Big(\frac{\partial\Pi}{\partial z}\Big)^2\Big]^{-1}
\Big(\frac{\partial\Pi}{\partial\varpi}\cdot\mbox{RH}\eqref{CEQd} +
\frac{\partial\Pi}{\partial z}\mbox{RH}\eqref{CEQe}\Big), \label{Z.K'a} \\
\frac{\partial K'}{\partial z}&=\Big[\Big(
\frac{\partial\Pi}{\partial\varpi}\Big)^2+
\Big(\frac{\partial\Pi}{\partial z}\Big)^2\Big]^{-1}
\Big(-\frac{\partial\Pi}{\partial z}\cdot\mbox{RH}\eqref{CEQd} +
\frac{\partial\Pi}{\partial \varpi}\mbox{RH}\eqref{CEQe}\Big), \label{Z.K'b}
\end{align}
\end{subequations}
where RH\eqref{CEQd}, RH\eqref{CEQe} stand for the right-hand sides of
\eqref{CEQd},\eqref{CEQe}, respectively, provided the assumption {\bf (B2)}, 
we can claim the following

\begin{Proposition}(\cite[Proposition 5]{asEE}) \label{Prop9}
Suppose {\bf (B2)}, and that $\Omega$ is a constant. Let $K'$ be arbitrarily fixed and let $F', A', \Pi, \rho$ satisfy
\eqref{CEQa} \eqref{CEQb} \eqref{CEQc} and \eqref{CEQf} with this fixed $K'$. Let us denote by
$\tilde{K}'_1, \tilde{K}'_3$ the right-hand sides of
\eqref{Z.K'a}, \eqref{Z.K'b}, respectively, evaluated by these $F', A', \Pi$. Then it holds that
\begin{align}
\frac{\partial\tilde{K}'_1}{\partial z}-
\frac{\partial\tilde{K}'_3}{\partial\varpi}&=
\frac{16\pi \mathsf{G}}{\mathsf{c}^4}
e^{2(-F'+K')}P\Pi\Big[\Big(\frac{\partial\Pi}{\partial\varpi}\Big)^2+
\Big(\frac{\partial\Pi}{\partial z}\Big)^2\Big]^{-1}\times \nonumber \\
&\times \Big[\Big(\frac{\partial K'}{\partial\varpi}-\tilde{K}'_1\Big)
\frac{\partial\Pi}{\partial z}-
\Big(\frac{\partial K'}{\partial z}-\tilde{K}'_3\Big)
\frac{\partial\Pi}{\partial\varpi}\Big].\label{3.39}
\end{align}
\end{Proposition}

Proof.  By a tedious calculation, we get
\begin{align*}
\frac{\partial\tilde{K}'_1}{\partial z}-\frac{\partial\tilde{K}'_3}{\partial\varpi}&=
-(\Pi_3\tilde{K}'_1-\Pi_1\tilde{K}'_3)
(\Pi_1^2+\Pi_3^2)^{-1}\cdot \triangle\Pi + \\
&+(\Pi_1^2+\Pi_3^2)^{-1}Z'
\end{align*}
with
\begin{align*}
Z':=&\frac{1}{2}(-\Pi_1\partial_3\triangle\Pi+\Pi_3\partial_1\triangle\Pi)+ \\
&+2\Pi(-\Pi_1F_3'+\Pi_3F_1')[S'\mathrm{a}]
+\frac{e^{4F'}}{2\Pi}(\Pi_1A_3'-\Pi_3A_1')[S'\mathrm{b}],
\end{align*}
where $\Pi_j, F_j', A_j', (j=1,3,)  \triangle\Pi$ stand for $\partial_j\Pi, \partial_jF', \partial_jA', \partial_1^2\Pi+\partial_3^2\Pi$ respectively, and $[S'\mathrm{a}], [S'\mathrm{b}](=0)$ stand for the right-hand sides of the equations \eqref{CEQa}, \eqref{CEQb}. 
Keeping in mind that
$$\partial_j\triangle\Pi=\Big[
2(-\partial_jF'+\partial_jK')-(\epsilon+P)\frac{\partial_jF'}{P}+\frac{\partial_j\Pi}{\Pi}\Big]\triangle\Pi,
$$
we get
$$Z'=(-\Pi_1\partial_3K'+\Pi_3\partial_1K')\triangle\Pi
$$
and therefore
\begin{align*}
\frac{\partial\tilde{K}'_1}{\partial z}-\frac{\partial\tilde{K}'_3}{\partial\varpi}&=
\Big[\Pi_1(\tilde{K}'_3-\partial_3K')-\Pi_3(\tilde{K}'_1-
\partial_1K')\Big]\times \\
&\times (\Pi_1^2+\Pi_3^2)^{-1}\triangle \Pi.
\end{align*}
This completes the proof. $\square$\\

Using this  Proposition \ref{Prop9}, we can claim the following theorem. The proof is the same as that of Theorem \ref{ThN2}.

\begin{Theorem} (\cite[Lemma 1]{asEE})
Let us consider a bounded bounded domain  $\mathfrak{D}$  and denote by $\bar{\mathfrak{D}}$ the closure of $\mathfrak{D}$. Assume that $\Omega$ is a constant on $\mathfrak{D}$. Suppose that $K' \in C^1(\bar{\mathfrak{D}})$ is given and that $F', A', \Pi, \rho \in C^3(\bar{\mathfrak{D}})$ satisfy \eqref{CEQa},\eqref{CEQb},\eqref{CEQc} and \eqref{CEQf} with this $K'$. 
Suppose that the assumption  {\bf (B2)} holds on $\bar{\mathfrak{D}}$. Let us denote by
$\tilde{K}'_1,\tilde{K}'_3$ the right-hand sides of
\eqref{Z.K'a},\eqref{Z.K'b}, respectively, evaluated by these
$F', A', \Pi$. (They are $C^1$-functions on $\bar{\mathfrak{D}}$.)
Put
\begin{equation}
\tilde{K}'(\varpi, z):=K_O'+
\int_0^z\tilde{K}_3'(0,z')dz'+
\int_0^{\varpi}\tilde{K}'_1(\varpi',z)d\varpi'\label{3.40}
\end{equation}
for $(\varpi, z) \in \mathfrak{D}$. Here $K_O'$ is a constant. If $\tilde{K}'=K'$, then $K'$
satisfies 
\begin{equation}
\frac{\partial K'}{\partial\varpi}=\tilde{K}'_1,\qquad
\frac{\partial K'}{\partial z}=\tilde{K}'_3,\label{3.41}
\end{equation}
that is, the equations \eqref{CEQd}\eqref{CEQe} are satisfied.
\end{Theorem}

\begin{Remark}
When $\Omega$ is not a constant in the considered domain, it seems that \eqref{3.39}
does not hold because of the presence of derivatives of $\Omega$.
\end{Remark}

\section{Corotating coordinate system and the meaning of $F', A', \Pi, K'$}

 In this section we assume that $\Omega$ is a constant on the domain to be considered.

We can consider the `{\bf corotating coordinate system}' characterized by
\begin{equation}
t'=t,\quad \varpi'=\varpi,\quad
\phi'=\phi-\Omega t,\quad z'=z. 
\end{equation}

It can be verified that the line element retains its form \eqref{2}  with
\begin{subequations}
\begin{align}
& e^{2F'}=e^{2F}\Big(1+A\frac{\Omega}{\mathsf{c}}\Big)^2
-e^{-2F}\Pi^2\frac{\Omega^2}{\mathsf{c}^2}, \label{X1}
\\
\mbox{or equivalently,}& \nonumber \\
&e^{2F}=e^{2F'}\Big(1-A'\frac{\Omega}{\mathsf{c}}\Big)^2
-e^{-2F'}\Pi^2\frac{\Omega^2}{\mathsf{c}^2},  \nonumber \\
\mbox{and} & \nonumber \\
& \Big(1-A'\frac{\Omega}{\mathsf{c}}\Big)e^{2F'}=
\Big(1+A\frac{\Omega}{\mathsf{c}}\Big)e^{2F}, \label{X6}\\
& K'-F'=K-F, \label{X4}
\\ 
&\Pi'=\Pi \label{X5}
\end{align}
\end{subequations}\\

Just to make sure, let us verify this observation.

First we define $F', A'$ by \eqref{X1}
and
\begin{equation}
e^{2F'}A'=e^{2F}\Big(1+\frac{\Omega}{\mathsf{c}}A\Big)A-
e^{-2F}\frac{\Omega}{\mathsf{c}}\Pi^2, \label{X2}
\end{equation}
provided that the right-hand side of \eqref{X1} is positive.
Then calculating
$$e^{-4F}\times\Big(\mbox{the square of \eqref{X2} }-\Pi^2\Big),
$$
we see that 
\begin{equation}
e^{2F'}(A')^2-e^{-2F'}\Pi^2=
e^{2F}A^2-e^{-2F}\Pi^2 \label{X3}
\end{equation}
holds. Therefore we can claim
\begin{equation}
ds^2=e^{2F'}
(\mathsf{c}dt'+A'd\phi')^2-
e^{-2F'}[e^{2K'}
((d\varpi')^2+(dz')^2)+(\Pi')^2(d\phi')^2], 
\end{equation}
by defining $K', \Pi'$ by \eqref{X4}, \eqref{X5}.
Moreover, calculating 
$$ \mbox{\eqref{X1}} -\frac{\Omega}{\mathsf{c}}\times\mbox{\eqref{X2}}, $$
we see that \eqref{X6} holds.\\

The primed quantities $\mathsf{c}^2F', \mathsf{c}A'$, etc are called `{\bf corotating potentials}'.\\

\begin{Remark}
Note that $F', A', K' ,  \Pi'(=\Pi)$ are the same as the  quantities defined in Section 4 by \eqref{3.21}, but, when $\Omega$, a function of $(\varpi, z)$, is not a constant, the metric generally does not retain the Lanczos form for the coordinate system $(t', \varpi', \phi', z')=(t, \varpi, \phi-\Omega t, z)$.
In fact we have
\begin{align*}
ds^2&=e^{2F}(\mathsf{c}t+Ad\phi)^2-e^{-2F}
(e^{2K}(d\varpi^2+dz^2)+\Pi^2d\phi^2) \\
&=e^{2F'}(\mathsf{c}t+A'd\phi')^2-e^{-2F'}
(e^{2K'}(d\varpi^2+dz^2)+\Pi^2(d\phi')^2) +\\
&+(e^{2F'}(A')^2-e^{-2F'}\Pi^2)
\sum_{j=1,3}\Big[
t^2(\partial_j\Omega)^2(dx^j)^2+
2t\Omega(\partial_j\Omega)dtdx^j+
2t(\partial_j\Omega)d\phi'dx^j\Big],
\end{align*}
where $\partial_j\Omega =\partial\Omega/\partial x^j, j=1,3$ and $x^1=\varpi, x^3=z$.
\end{Remark}

Since 
\begin{equation}
e^{2F}\Big(1+A\frac{\Omega}{\mathsf{c}}\Big)^2
-e^{-2F}\Pi^2\frac{\Omega^2}{\mathsf{c}^2}=e^{2G}, 
\end{equation}
 looking at \eqref{X1}, we see that $G$ is equal to the corotating
potential $F'$, that is,
\begin{equation}
G=F'.
\end{equation}\\

The components of the 4-velocity $U^{\mu'}$ with respect to the corotating coordinate system
$$x^{0'}=\mathsf{c} t',\quad x^{1'}=\varpi',\quad
x^{2'}=\phi',\quad x^{3'}=z' $$
turn out to be
\begin{equation}
U^{0'}=e^{-G},\quad U^{1'}=U^{2'}=U^{3'}=0,
\end{equation}
or
\begin{equation}
U^{\mu'}\frac{\partial }{\partial x^{\mu'}}=e^{-G}\frac{1}{\mathsf{c}}\frac{\partial}{\partial t'}.
\end{equation}

\vspace{20mm}
{\bf\Large Acknowledgment}\\

This work was partially done on the occasion of the BIRS-CMO Workshop `Time-like Boundaries in General Relativistic Evolution Problems (19w5140)' held at Oaxaca, Mexico on July 28 - August 2, 2019. The author would like to express his thanks to Professors Helmut Friedrich, Olivier Sarbach, Oscar Reula for organizing this workshop, to the Casa Matematica Oaxaca of the Banff International Research Station for Mathematical Innovation and Discovery for the hospitality, and to the participants for discussions. This work is supported by JPS KAKENHI Grant Number JP18K03371.

\vspace{20mm}


{\bf\Large  Appendix 1}\\

The components of the metric and the Christoffel symbols 
$$ \Gamma_{\nu\lambda}^{\mu}=\frac{1}{2}g^{\mu\alpha}
(\partial_{\lambda}g_{\alpha\nu}+
\partial_{\nu}g_{\alpha\lambda}-\partial_{\alpha}g_{\nu\lambda} )$$
are
as following ( other $g_{\mu\nu}, g^{\mu\nu}, \Gamma_{\mu\nu}^{\lambda}$ are zero ):

\begin{align*}
&g_{00}=e^{2F}, \quad g_{02}=g_{20}=e^{2F}A, \\
&g_{11}=g_{33}=-e^{-2F+2K}, \quad g_{22}=e^{2F}A^2-e^{-2F}\Pi^2;
\end{align*}
\begin{align*}
& g^{00}=-\frac{1}{\Pi^2}(e^{2F}A^2-e^{-2F}\Pi^2), \quad g^{02}=g^{20}=\frac{e^{2F}}{\Pi^2}A, \\
& g^{11}=g^{33}=-e^{2F-2K}, \quad g^{22}=-\frac{e^{2F}}{\Pi^2};
\end{align*}
\begin{align*}
 \Gamma_{01}^0&=\Gamma_{10}^0=\frac{1}{2\Pi^2}
(e^{4F}A\partial_1A+\Pi^2e^{-2F}\partial_1e^{2F}) =\frac{e^{4F}}{2\Pi^2}A\partial_1A+\partial_1F, \\
 \Gamma_{03}^0&=\Gamma_{30}^0=\frac{1}{2\Pi^2}
(e^{4F}A\partial_3A+\Pi^2e^{-2F}\partial_3e^{2F}) =\frac{e^{4F}}{2\Pi^2}A\partial_3A+\partial_3F, \\
 \Gamma_{12}^0&=\Gamma_{21}^0=
\frac{1}{2\Pi^2}(\Pi^2\partial_1A-
e^{4F}A\partial_1(e^{-4F}\Pi^2)) =\frac{1}{2}\partial_1A-\frac{A}{\Pi}\partial_1\Pi+2A\partial_1F, \\
\Gamma_{23}^0&=\Gamma_{32}^0=
\frac{1}{2\Pi^2}(\Pi^2\partial_3A-
e^{4F}A\partial_3(e^{-4F}\Pi^2)) =\frac{1}{2}\partial_3A-\frac{A}{\Pi}\partial_3\Pi+2A\partial_3F
\end{align*}
\begin{align*}
\Gamma_{00}^1&=\frac{1}{2}e^{2F-2K}\partial_1e^{2F}=e^{4F-2K}\partial_1F, \\
\Gamma_{02}^1&=\Gamma_{20}^1=\frac{1}{2}e^{2F-2K}
\partial_1(e^{2F}A),\\
\Gamma_{11}^1&=\frac{1}{2}e^{2F-2K}\partial_1e^{-2F+2K}
=-\partial_1(F-K), \\
\Gamma_{13}^1&=\Gamma_{31}^1=
\frac{1}{2}e^{2F-2K}\partial_3e^{-2F+2K}
=-\partial_3(F-K), \\
\Gamma_{22}^1&=\frac{1}{2}e^{2F-2K}
\partial_1(e^{2F}A^2-e^{-2F}\Pi^2), \\
\Gamma_{33}^1&=-\frac{1}{2}
e^{2F-2K}\partial_1e^{-2F+2K}
=\partial_1(F-K);
\end{align*}
\begin{align*}
\Gamma_{01}^2&=\Gamma_{10}^2=-\frac{e^{4F}}{2\Pi^2}
\partial_1A, \\
\Gamma_{03}^2&=\Gamma_{30}^2=-\frac{e^{4F}}{2\Pi^2}
\partial_3A, \\
\Gamma_{12}^2&=\Gamma_{21}^2=
\frac{e^{2F}}{2\Pi^2}
(-e^{2F}A\partial_1A+\partial_1(e^{-2F}\Pi^2)) 
=-\frac{e^{4F}}{2\Pi^2}A\partial_1A-\partial_1F+\frac{1}{\Pi}\partial_1\Pi, \\
\Gamma_{23}^2&=\Gamma_{32}^2=
\frac{e^{2F}}{2\Pi^2}
(-e^{2F}A\partial_3A+\partial_3(e^{-2F}\Pi^2)) 
=-\frac{e^{4F}}{2\Pi^2}A\partial_3A-\partial_3F+\frac{1}{\Pi}\partial_3\Pi ;
\end{align*}
\begin{align*}
\Gamma_{00}^3&=\frac{1}{2}e^{2F-2K}\partial_3e^{2F}
=e^{4F-2K}\partial_3F, \\
\Gamma_{02}^3&=\Gamma_{20}^3=\frac{1}{2}e^{2F-2K}
\partial_3(e^{2F}A)
, \\
\Gamma_{11}^3&=-\frac{1}{2}e^{2F-2K}\partial_3e^{-2F+2K}
=\partial_3(F-K), \\
\Gamma_{13}^3&=\Gamma_{31}^3=
\frac{1}{2}e^{2F-2K}\partial_1e^{-2F+2K}
=-\partial_1(F-K), \\
\Gamma_{22}^3&=\frac{1}{2}e^{2F-2K}
\partial_3(e^{2F}A^2-e^{-2F}\Pi^2), \\
\Gamma_{33}^3&=\frac{1}{2}
e^{2F-2K}\partial_3e^{-2F+2K}
=-\partial_3(F-K).
\end{align*}\\

{\bf\Large Appendix 2}\\

The the components of the metric and the Christoffel symbols describing by $f, k, l, m$ are given 
as following (other
$g_{\mu\nu}, g^{\mu\nu}, \Gamma_{\mu\nu}^{\lambda}$ are zero ):
\begin{align*}
&g_{00}=f, \quad g_{02}=g_{20}=-k, \\
&g_{22}=-l, \quad g_{11}=g_{33}=-e^m;
\end{align*}
\begin{align*}
&g^{00}=\frac{l}{\Pi^2},\quad g^{02}=g^{20}=-\frac{k}{\Pi^2}, \\
&g^{22}=-\frac{f}{\Pi^2}, \quad g^{11}=g^{33}=-e^{-m};
\end{align*}
\begin{align*}
&\Gamma_{01}^0=\Gamma_{10}^0=\frac{1}{2\Pi^2}(l\partial_1f+k\partial_1k), \quad
\Gamma_{03}^0=\Gamma_{30}^0=\frac{1}{2\Pi^2}(l\partial_3f+k\partial_3k), \\
&\Gamma_{12}^0=\Gamma_{21}^0=\frac{1}{\Pi^2}(k\partial_1l-l\partial_1k), \quad
\Gamma_{23}^0=\Gamma_{32}^0=\frac{1}{2\Pi^2}(k\partial_3l-l\partial_3k);
\end{align*}
\begin{align*}
&\Gamma_{00}^1=\frac{1}{2}e^{-m}\partial_1f, \quad
\Gamma_{02}^1=\Gamma_{20}^1=-\frac{1}{2}e^{-m}\partial_1k, \quad
\Gamma_{11}^1=\frac{1}{2}\partial_1m, \\
&\Gamma_{13}^1=\Gamma_{31}^1=\frac{1}{2}\partial_3m, \quad
\Gamma_{33}^1=-\frac{1}{2}\partial_1m, \quad
\Gamma_{22}^1=-\frac{1}{2}e^{-m}\partial_1l;
\end{align*}
\begin{align*}
&\Gamma_{01}^2=\Gamma_{10}^2=\frac{1}{2\Pi^2}
(f\partial_1k-k\partial_1f), \quad
\Gamma_{03}^2=\Gamma_{30}^2=\frac{1}{2\Pi^2}
(f\partial_3k-k\partial_3f), \\
&\Gamma_{12}^2=\Gamma_{21}^2=\frac{1}{2\Pi^2}
(f\partial_1l+k\partial_1k), \quad
\Gamma_{23}^2=\Gamma_{32}^2=\frac{1}{2\Pi^2}
(f\partial_3l+k\partial_3k);
\end{align*}
\begin{align*}
&\Gamma_{00}^3=\frac{1}{2}e^{-m}\partial_3f, \quad
\Gamma_{02}^3=\Gamma_{20}^3=-\frac{1}{2}e^{-m}\partial_3k, \quad
\Gamma_{11}^3=-\frac{1}{2}\partial_3m, \\
&\Gamma_{13}^3=\Gamma_{31}^3=\frac{1}{2}\partial_1m, \quad
\Gamma_{33}^3=\frac{1}{2}\partial_3m, \quad
\Gamma_{22}^3=-\frac{1}{2}e^{-m}\partial_3l.
\end{align*}

\end{document}